\documentclass{article}

\usepackage[english]{babel}

\usepackage[letterpaper,top=2cm,bottom=2cm,left=3cm,right=3cm,marginparwidth=1.75cm]{geometry}

\usepackage{amsmath}
\usepackage{graphicx}
\usepackage[colorlinks=true, allcolors=blue]{hyperref}
\usepackage{subfigure}
\usepackage{caption}
\usepackage{lipsum}
\usepackage{booktabs}
\date{}
\begin{document}

\title{Fast Fourier Transform periodic interpolation method for superposition sums in a periodic unit cell}

\author{Fangzhou Ai \and Vitaliy Lomakin\thanks{Corresponding Author. Email: vlomakin@ucsd.edu}}

\maketitle

\begin{center}
    Department of Electrical and Computer Engineering, University of California San Diego, \\
    San Diego 92093, U.S.A \\
    Center for Memory and Recording Research, University of California San Diego, \\
    San Diego 92093, U.S.A
\end{center}

\begin{abstract}
We propose a Fast Fourier Transform based Periodic Interpolation Method (FFT-PIM), a flexible and computationally efficient approach for computing the scalar potential given by a superposition sum in a unit cell of an infinitely periodic array. Under the same umbrella, FFT-PIM allows computing the potential for 1D, 2D, and 3D periodicities for dynamic problems involving the Helmholtz potential and static problems involving Coulomb potential, including problems with and without a periodic phase shift. The computational complexity of the FFT-PIM is of $O(N \log N)$ for $N$ spatially coinciding sources and observer points. The FFT-PIM uses rapidly converging series representations of the Green's function serving as a kernel in the superposition sum. Based on these representations, the FFT-PIM splits the potential into its near-zone component, which includes a small number of images surrounding the unit cell of interest, and far-zone component, which includes the rest of an infinite number of images. The far-zone component is evaluated by projecting the non-uniform sources onto a sparse uniform grid, performing superposition sums on this sparse grid, and interpolating the potential from the uniform grid to the non-uniform observation points. The near-zone component is evaluated using an FFT-based method, which is adapted to efficiently handle non-uniform source-observer distributions within the periodic unit cell. The FFT-PIM can be used for a broad range of applications, such as periodic problems involving integral equations for wave propagation in electromagnetics and acoustics, micromagnetic solvers, and density functional theory solvers. 
\end{abstract}

\section{Introduction} \label{sec:1}

\begin{figure}[h] 
\includegraphics[width=8cm]{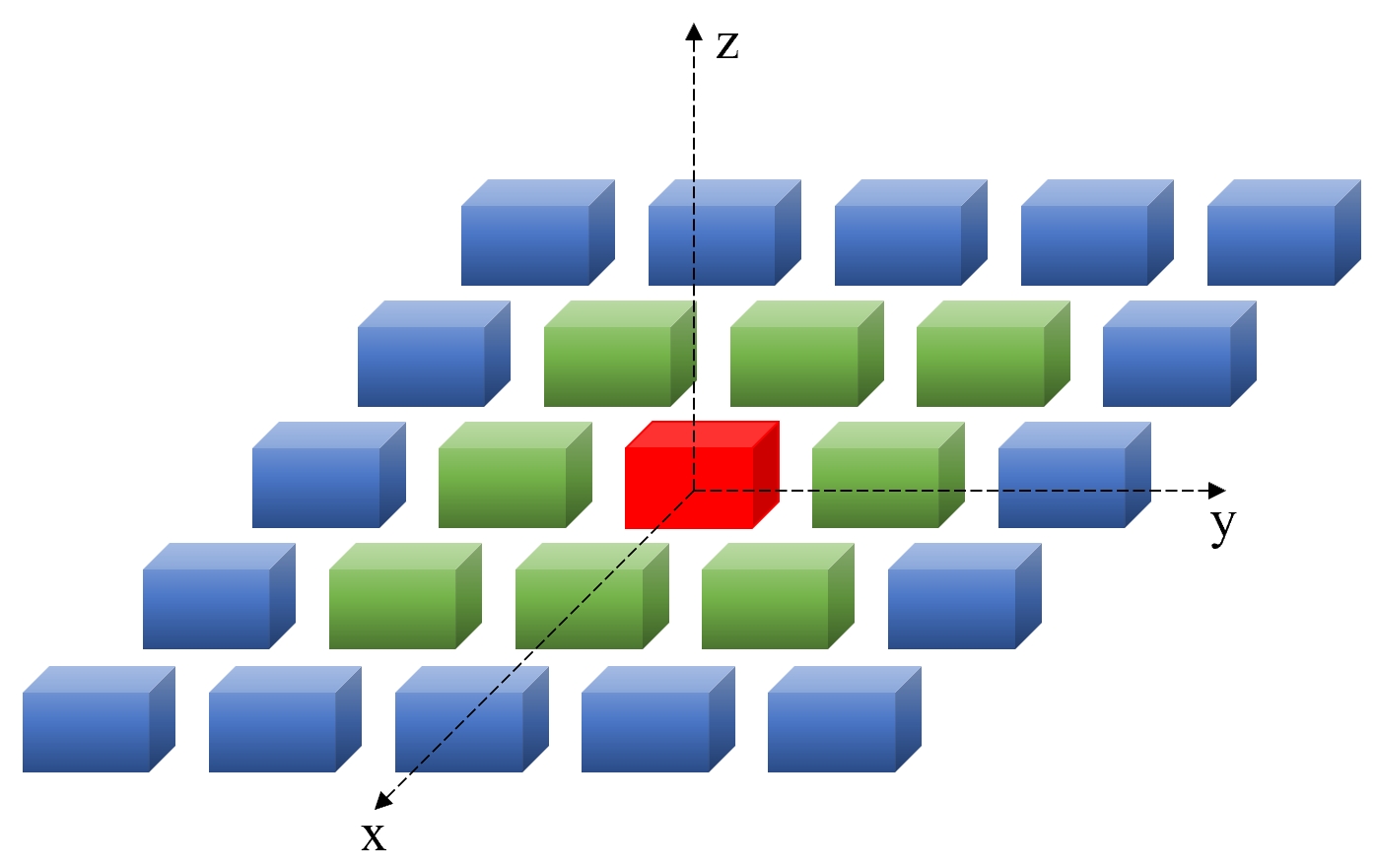}
\centering
\caption{Illustration for a 2D periodic problem consisting of an infinite 2D array of cubes. The central red cube is the zeroth unit cell. Surrounding it along the $x$ and $y$ axes are the $1^{st}$ order image cubes in green and the $2^{nd}$ order image cubes in blue.} \label{fig:illu_prob}
\end{figure}

A common category of problem in computational physics is calculating the potential, referred to as periodic scalar potential (PSP), that is generates by a non-uniform distribution of sources arranged in a unit cell that is extended to an infinite periodic array (Fig.~\ref{fig:illu_prob}). Examples of such problems are phased-array antennas, crystals, periodic gratings, periodically modulated waveguides, to name a few \cite{mailloux2017phased, kalkstein1971green, sirenko2010modern, peng1975theory}. These problems may present challenges for obtaining accurate and rapid solutions. For a \textcolor{black}{non-periodic} unit cell, this task can be accomplished via a standard superposition sum involving free-space Green's function as its kernel. Such sums can be evaluated rapidly in $O(N)$ or $O(N\log N)$ computational cost for $N$ sources and observers using several fast methods, such as fast multipole like method \cite{GREENGARD1987325, 6230628}, H-matrix method \cite{hackbusch2000sparse}, interpolation-based methods \cite{1159856, LI20108463, MENG20108430}, and FFT-based methods \cite{NUFFT, ali_yilmaz, 7770011, 662670}. However, for the periodic case, the superposition sums involve infinite number of images, which poses additional challenges to the slow converge or even divergence as well as a potentially high computational cost. 

Several existing methods can address such periodic problems with certain assumptions and limitations\textcolor{black}{\cite{r1_1, r1_2, r1_3, r1_4}}. Periodic Green's function (PGF) can be defined in terms of an infinite sum over periodic images. In dynamic problems, e.g., in computational electromagnetics, PGF is assumed to be defined with a phase shift or non-vanishing wavenumber \cite{CAPOLINO2007250, FPIP, Derek}. However, for static problems with no phase shift, even in the 1D case, the PGF diverges. In molecular dynamics problems \cite{molecule}, a special form of static PGF assumes uniformly arranged sources \cite{LGF}, such as a cubic lattice. Although various methods have been proposed to address static and dynamic periodic unit cell problems with non-uniform sources, they often focus on specific problem types, such as addressing only 1D or 2D periodicity, addressing only static problems \cite{Lebecki_2008}, requiring lattice structure or relying on finite difference methods \cite{fdtd}. There is a need in a general approach capable of effectively handling all types of these problems regardless of the periodicity type, distribution of sources, with or without phase shift, and for static and dynamic cases.

In this paper, we proposed a general approach, referred to as Fast Fourier Transform Periodic Interpolation Method (FFT-PIM) that is applicable to all these types of problems under a unified framework, including 1D, 2D, and 3D periodicities for static and dynamic problems with and without periodic phase shifts. The FFT-PIM separates the PSP into a near-zone component involving a small number of images near the periodic unit cell of interest, and a spatially slowly varying far-zone component involving the rest of an infinite number of images. The evaluation of the near-zone PSP component is based on the box adaptive integral method \cite{NUFFT}, which is modified to allow for the rapid computations with multiple near-zone images. The evaluation of the far-zone PSP component is based on the sparse periodic interpolation method \cite{5582244}, which is adapted to allow handling dynamic and static problems with arbitrary phase shifts. The computational cost of the FFT-PIM is of $O(N\log N)$ and the memory consumption is of $O(N)$.

The paper is organized as follows. Sec.~\ref{sec:2} presents the problem formulation. Sec.~\ref{sec:3} presents the algorithmic foundations of the FFT-PIM, including near- and far-zone representations of PGF and PSP as well as efficient ways of their evaluation. Sec.~\ref{sec:4} shows results for the use of the FFT-PIM, its error analysis, and computational performance. Finally, Sec.~\ref{sec:5} presents summary and conclusions. 

\section{Formulation} \label{sec:2}
Consider an infinite periodic array of unit cells residing in free space (see Fig.~\ref{fig:illu_prob} for a 2D periodicty example). Within each unit cell there are $N$ source points located at $\mathbf{r}_n$. In the zeroth unit cell, the source amplitudes are $q_n$, where $n=1...N$. In each unit cell, there is the same number of coinciding observer points at which the periodic scalar potentials (PSPs) $u(\mathbf{r}_i)$ are to be found. The periodic array can be 1D, 2D or 3D with the periodicity of $L_x$, $L_y$, and $L_z$ in three possible directions $x$, $y$, and $z$, respectively. The amplitudes of the sources may be periodically phase shifted with a linear, possibly complex, phase shift determined by the wavenumbers $k_{x0}$, $k_{y0}$, and $k_{z0}$ along the $x$, $y$, and $z$ directions. A free space wavenumber $k_0$ describes the propagation of the waves in the free space for the dynamic case. For $k_0=0$, the problem is static, and for $k_{x0}=k_{y0}=k_{z0}=0$, the sources are not phase shifted.

The PSP in the zeroth unit cell can be written as 
\begin{equation}
    \label{eq:1}
    u(\mathbf{r}_m)=\sum_{n=1}^N G^p(\mathbf{r}_m-\mathbf{r}_n)q(\mathbf{r}_n),
\end{equation}
where $m=1...N$ and $G^p$ represents the scalar periodic Green's function (PGF), given for the 1D, 2D, and 3D cases as:
\begin{subequations}
    \label{eq:2}
    \begin{align}
        G^{p}_{1D}(\mathbf{r}) =& \sum_{i_x=-\infty}^{\infty} e^{-jk_{x0}(i_xL_x)}G_0(\mathbf{r}-i_xL_x\mathbf{\hat{x}}),\\
        G^{p}_{2D}(\mathbf{r}) =& \sum_{i_x=-\infty}^{\infty}\sum_{i_y=-\infty}^{\infty}e^{-j[k_{x0}(i_xL_x)+k_{y0}(i_yL_y)]}G_0(\mathbf{r}-i_xL_x\mathbf{\hat{x}}-i_yL_y\mathbf{\hat{y}}),\\
        G^{p}_{3D}(\mathbf{r}) =& \sum_{i_x=-\infty}^{\infty}\sum_{i_y=-\infty}^{\infty}\sum_{i_z=-\infty}^{\infty}e^{-j[k_{x0}(i_xL_x)+k_{y0}(i_yL_y)+k_{z0}(i_zL_z)]}G_0(\mathbf{r}-i_xL_x\mathbf{\hat{x}}-i_yL_y\mathbf{\hat{y}}-i_zL_z\mathbf{\hat{z}}),
    \end{align}
\end{subequations}
respectively. Here, $G_0$ is the free space scalar Green's function given by $G_0(\mathbf{r})=\exp(-jk_0|\mathbf{r}|)/4\pi|\mathbf{r}|$. 

The series in Eq.~\eqref{eq:1} for the PSP is important in multiple computational physics problems. For example, the evaluation of the superposition integrals appearing in electromagnetic integral equation or micromagnetic solvers \cite{Fastmag} can be represented as a product of sparse matrices describing the mesh of the geometry and the result of the superposition sum of Eq.~\eqref{eq:1}. The sum of Eq.~\eqref{eq:1} can also be directly apply to a discrete set of sources and observers, e.g., when considering a set of point sources. 

The series of Eq.~\eqref{eq:2} for the PGF is slowly convergent or even divergent depending on the wavenumbers, which is because of a slow spatial decay of $G_0$. An alternative spectral representation of the PGF of Eq.~\eqref{eq:2} can be based on the Floquet mode expansion \cite{Derek}:
\begin{subequations}\label{eq:pgf}
    \begin{align}
        G^{p}_{1D}(\mathbf{r})=&\sum_{m=-\infty}^{\infty}\frac{1}{4jL_x}e^{-jk_{xm}x}H_0^{(2)}(k_{\rho m}\sqrt{y^2+z^2}), \label{eq:pgf_1} \\
        G^{p}_{2D}(\mathbf{r})=&\sum_{m,n=-\infty}^{\infty}\frac{e^{-jk_{xm}x-jk_{yn}y-jk_{zmn}|z|}}{2jk_{zmn}L_xL_y}, \label{eq:pgf_2} \\
        G^{p}_{3D}(\mathbf{r})=&\sum_{m,n=-\infty}^{\infty}\frac{e^{-jk_{xm}x-jk_{yn}y}}{2jk_{zmn}L_xL_y}\nonumber \\
        &\times\left(e^{-jk_{zmn}|z|}+\frac{e^{-j(k_{zmn}-k_{z0})L_z}e^{-jk_{zmn}z}}{1-e^{-j(k_{zmn}-k_{z0})L_z}}+  \frac{e^{-j(k_{zmn}+k_{z0})L_z}e^{jk_{zmn}z}}{1-e^{-j(k_{zmn}+k_{z0})L_z}}\right), \label{eq:pgf_3}
    \end{align}
\end{subequations}
where $H_0^{(2)}$ is the zeroth order Hankel function of the second kind, $k_{xm}=k_{x0}+2\pi m /L_x$, $k_{ym}=k_{y0}+2\pi n /L_z$, $k_{\rho m}=(k_0^2-k_{xm}^2)^{1/2}$ and $k_{zmn}=(k_0^2-k_{xm}^2-k_{yn}^2)^{1/2}$, and the square roots are chosen such that their imaginary parts are non-positive. The spectral series of Eq.~\eqref{eq:pgf} converges exponentially fast when $(y^2+z^2)^{1/2}$ is non-vanishing and $k_0,k_{x0},k_{y0},k_{z0}$ are not all vanishing. \textcolor{black}{The series may diverge or not be defined for certain parameter combinations with $k_{zmn}=0$, $k_{\rho m}=0$, which are related to so called Rayleigh-Wood anomalies \cite{Barnett2011, PhysRevB.71.235117, 4232642, Hessel:65}. For scattering electromagnetic or acoustic problems, such anomalies correspond to the transition between evanescent and propagating Floquet modes. Additional types of Wood anomalies may be of resonant type \cite{Hessel:65}, e.g., when Green's function is defines in the presence of a layered medium \cite{6059491}. When such anomalies occur, they may require a special treatment, e.g., Ref.~\cite{Barnett2011} or by defining the phase shift in the complex plane with the integration path deformation \cite{4232642,1603823}.}

A special attention needs to be paid to the \textcolor{black}{no phase static periodic (NPSP) case} where $k_0=k_{x0}=k_{y0}=k_{z0}=0$. This case is important as it corresponds to several practical problem, such as using periodic extensions to represent infinite 1D, 2D, and 3D domains in various static problems, e.g., to characterize micromagnetic \cite{10.1063/1.348954} or in molecular dynamics problems. In this \textcolor{black}{NPSP} case, the PGF calculated via series of Eq.~\eqref{eq:2} or Eq.~\eqref{eq:pgf} diverges and the PSP for a general source distribution diverges as well. However, a finite PSP can be obtained for the special case of a neutral source, i.e., under the condition
\begin{equation}
    \label{eq:neutral}
    \sum_{l=1}^N q_l = 0.
\end{equation}
This neutrality condition corresponds to various physical problems, e.g., when computing the magnetic scalar potential generated by a magnetized objects or electrostatic potential in polarized objects \cite{stohr2006magnetism}. In this case, the PSP converges and we can neglect any non-relevant components (e.g., constant terms) in the PGF, so that it is convergent as well. The result is that the PGF for the neutral \textcolor{black}{NPSP} case can be calculated via:
\begin{subequations} \label{eq:lgf}
    \begin{align}
         G^{p}_{1D}(\mathbf{r})=&-\frac{1}{2\pi L_x}\ln\left(\sqrt{y^2+z^2}\right)+\frac{1}{\pi L_x}\sum_{m=1}^{\infty}K_0\left(\frac{2m\pi}{L_x} \sqrt{y^2+z^2}\right) \cos \left(\frac{2m\pi x}{L_x}\right) \label{eq:lgf_1}\\
         G^{p}_{2D}(\mathbf{r})=&-\frac{|z|}{2L_xL_y} -\frac{1}{4\pi L_x}\ln\left(1-2e^{-\frac{2\pi|z|}{L_y}}\cos \left(\frac{2\pi y}{L_y}\right)+e^{-\frac{4\pi|z|}{L_y}}\right) \nonumber \\
         & +\frac{1}{\pi L_x}\sum_{m=1}^{\infty}\sum_{n=-\infty}^{\infty}K_0\left(\frac{2m\pi}{L_x}\sqrt{(nL_y+y)^2+z^2}\right) \cos \left(\frac{2m\pi x}{L_x}\right) \label{eq:lgf_2}\\
         G^{p}_{3D}(\mathbf{r})=&\frac{1}{2L_xL_yL_z}(z^2-|z|L_z) \nonumber \\
         &-\frac{1}{4\pi L_x}\sum_{k=-\infty}^{\infty}\ln\left(1-2e^{-\frac{2\pi|kL_z+z|}{L_y}}\cos \left(\frac{2\pi y}{L_y}\right)+e^{-\frac{4\pi|kL_z+z|}{L_y}}\right)\nonumber \\
         &+\frac{1}{\pi L_x}\sum_{m=1}^{\infty}\sum_{n,k=-\infty}^{\infty}K_0\left(\frac{2m\pi}{L_x}\sqrt{(nL_y+y)^2+(kL_z+z)^2}\right) \cos \left(\frac{2m\pi x}{L_x}\right) \label{eq:lgf_3}
    \end{align}    
\end{subequations}
where $K_0$ is the zeroth order modified Bessel function of the second kind. 

It can be shown that Eq.~\eqref{eq:lgf_1} is obtained from Eq.~\eqref{eq:pgf_1} under the assumption of Eq.~\eqref{eq:neutral}. In this case, by taking the limit of ($k_0\to0, k_{x0}\to 0$) for all $n\neq 0$, the Hankel function becomes the modified Bessel function and the symmetric $\pm n$ terms can be combined into the positive $n$ terms. For the $n=0$ term, the Hankel function exhibits a logarithmic behavior but it is summed up to $0$ due to the source neutrality of Eq.~\eqref{eq:neutral}. The PGF for the 2D and 3D cases can be then obtained by a spatial sum of 1D PGFs. When the sources are arranged on a lattice, the PGF for the \textcolor{black}{NPSP} case sometimes is also referred to as the lattice Green's function \cite{LGF}.

\begin{figure}[htbp]
    \subfigure[]{
    \includegraphics[width=0.47\linewidth]{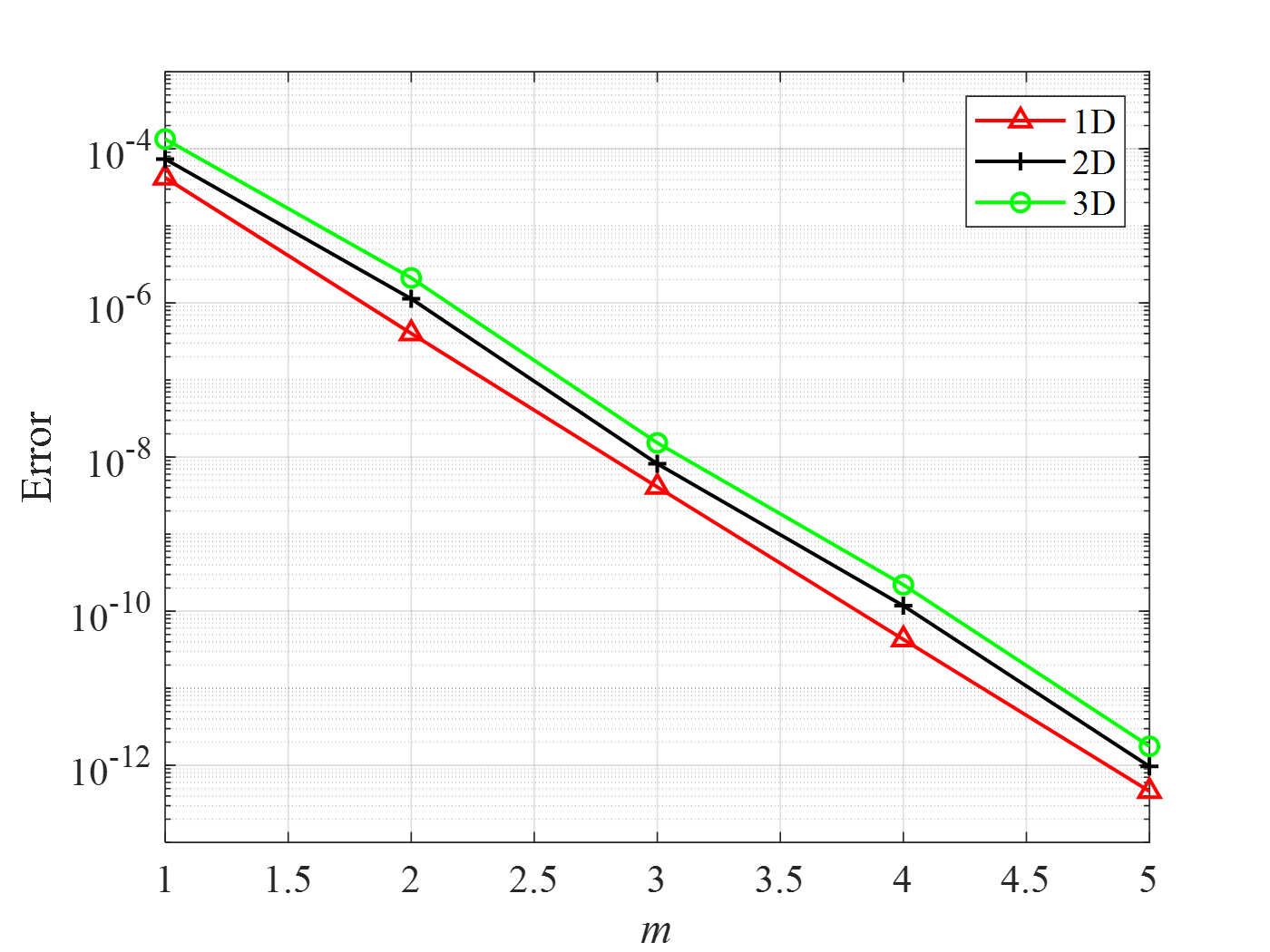}\label{fig:dpgf}
    }
    \quad
    \subfigure[]{
    \includegraphics[width=0.47\textwidth]{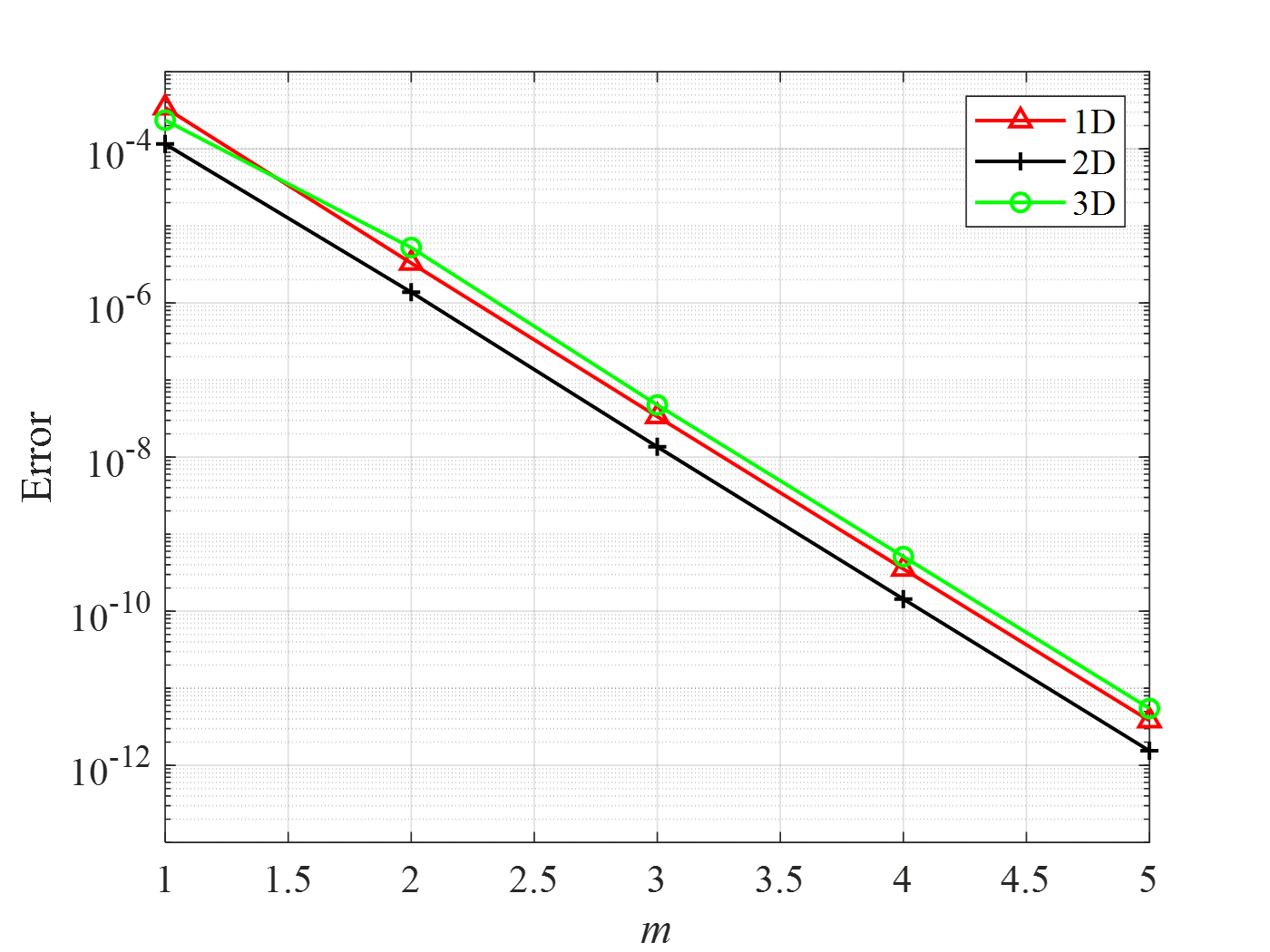}\label{fig:dlgf}
    }
    \caption{Convergence of PGFs, $L_x,L_y,L_z=1$, $x=y=z=0.5$. (a) The relative error of the sum of the first $m$ terms from Eq.~\eqref{eq:pgf}, with $k_0=-1-j$ and $k_{x0}=1-j, k_{y0}=1+j, k_{z0}=-1+j$, (b) The relative error of the sum of the first $m$ terms from Eq.~\eqref{eq:lgf}. Exponential convergence is achieved in both cases.}
    \label{fig:dpgf_dlgf}
\end{figure}

The PGFs converge exponentially fast provided the $(y^2+z^2)^{1/2}$ is not too small. For a desired error of $\epsilon$, the number of terms required Eq.~\eqref{eq:pgf_1} and Eq.~\eqref{eq:lgf_1} scales as $L_x\log(\epsilon^{-1})/z$ for the 1D case and $L_xL_y\log^2(\epsilon^{-1})/z^2$ for the 2D and 3D cases. To demonstrate the convergence, Fig.~\ref{fig:dpgf_dlgf} shows the PGFs and their convergence when using Eq.~\eqref{eq:pgf} and Eq.~\eqref{eq:lgf} for $L_x,L_y,L_z=1$, $k_0=-1-j$ and $k_{x0}=1-j, k_{y0}=1+j, k_{z0}=-1+j$ in Fig.~\ref{fig:dpgf} and the \textcolor{black}{NPSP} case in Fig.~\ref{fig:dlgf}.

\section{Algorithm} \label{sec:3}
\subsection{Representation of PGF and PSP} \label{ssec:3.1}
We now describe the FFT-PIM that allows rapidly computing the PSP in Eq.~\eqref{eq:1}. We start by representing the PGF in terms of its near-zone component $G_{near}^{p}$ and far-zone component $G^p_{far}$:
\begin{equation} \label{eq:4}
    G^p = G_{near}^{p} + G^p_{far}.
\end{equation}
Here, the near-zone component $G_{near}^{p}$ is given by Eq.~\eqref{eq:2} but with the sums over $i_x$, $i_y$, and $i_z$ in the range of $[-i_{xd},i_{xd}]$, $[-i_{yd},i_{yd}]$, and $[-i_{zd},i_{zd}]$, respectively, where $i_{xd},i_{yd},i_{zd}$ are small integer numbers. This components represent the contribution from a small number of images around the zeroth unit cell of interest. Assuming that the periodic lengths are close to each other, we can use the same value of $i_d$ for $i_{xd},i_{yd},i_{yd}$. For the sake of clarity we assume this case in the following discussions. $G_{near}^{p}$ includes the self-term and a few terms from the surrounding images (see Fig.~\ref{fig:illu_prob} showing the near-zone images for the case of $i_{xd}=i_{yd}=1$ as green boxes). Therefore, $G_{near}^{p}$ may be spatially singular and rapidly varying. The evaluation of $G_{near}^{p}$ involves only a small number of the sum terms. 

The far-zone component $G_{far}^{p}$ is the rest of the images from far-away unit cells and can be evaluated as $G^p_{far} = G^p - G_{near}^{p}$ (blue boxes in Fig.~\ref{fig:illu_prob}). The far-zone PGF $G_{far}^{p}$ is smoothly varying in space because it corresponds to large source-observer spatial separations. The rate of the variations is reduced by using a greater $i_d$.

Following the near- and far-zone decomposition of the PGF, the potential is also expressed in terms of its near- and far-zone components:
\begin{equation} \label{eq:5}
    \begin{split}
        u(\mathbf{r}_m)& = u_{near}(\mathbf{r}_m)+u_{far}(\mathbf{r}_m),\\
        u_{near}(\mathbf{r}_m)&=\sum_{n=1,n\neq m}^NG_{near}^{p}(\mathbf{r}_m- \mathbf{r}_n)q_n,\\
        u_{far}(\mathbf{r}_m)&=\sum_{n=1,n\neq m}^NG^p_{far}(\mathbf{r}_m- \mathbf{r}_n)q_n,
    \end{split}
\end{equation}
where, again, the far-zone potential $u_{far}$ is a spatially slowly varying function. These potential and PGF decompositions lead to a fast approach for the numerical evaluation of the PSP $u(\mathbf{r}_m)$.

\subsection{Evaluation of the far-zone PSP} \label{ssec:3_far_field}
\begin{figure}[h] 
\includegraphics[width=8cm]{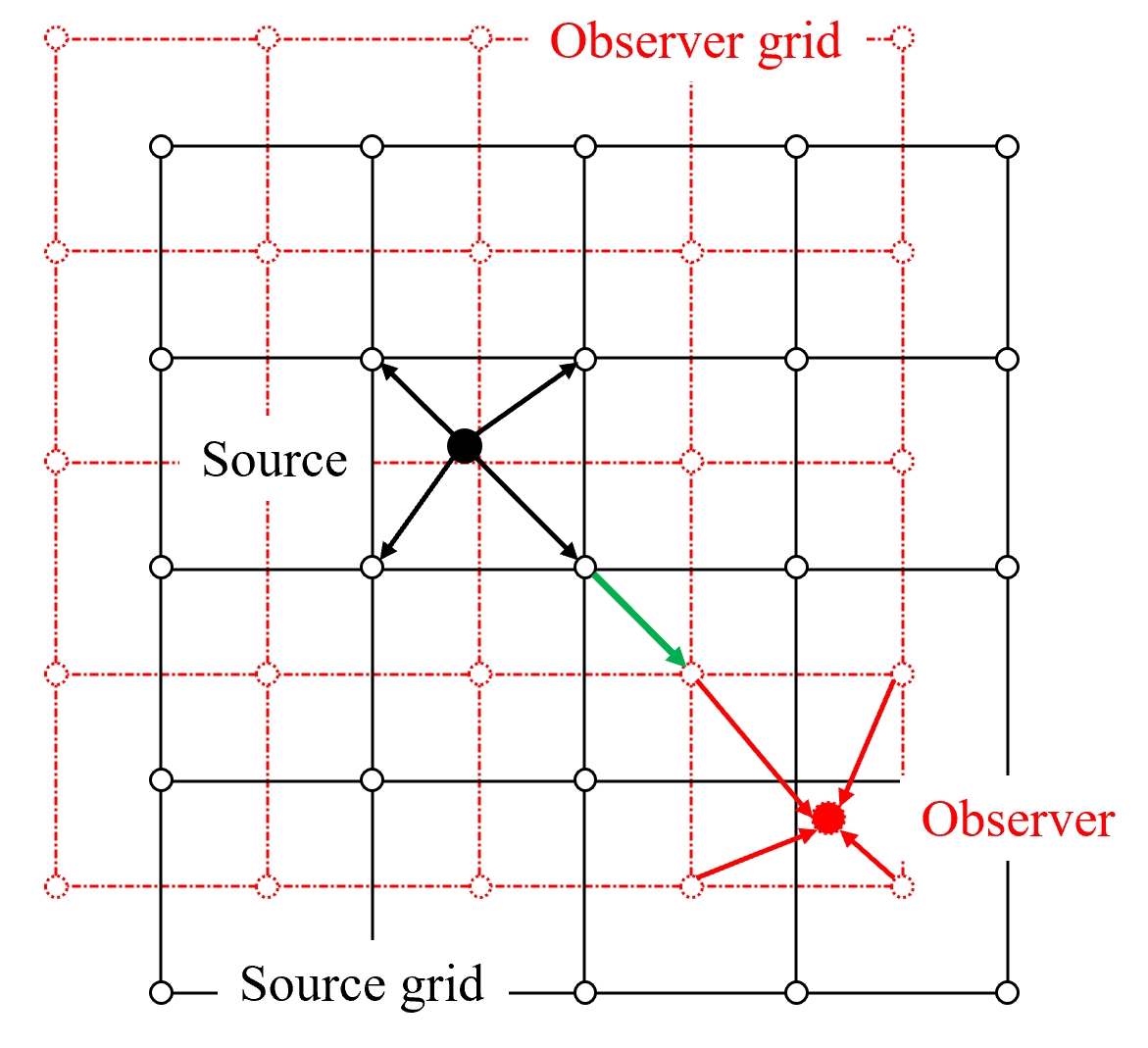}
\centering
\caption{Illustration of the source (black) and observer (red) grids for the far-zone PSP component. Black and red circles are source/observer grid points on the Cartesian lattices. Black and red dots are source/observer points with arbitrary distribution. Black arrows denote projections from source points to source grid points. Green arrow is the grid-to-grid interaction, performed by convolution. Red arrows represent interpolations from observer grid points to observer points.} \label{fig:grid}
\end{figure}

The far-zone PGF and PSP components vary slowly in space, which enables us to calculate the far-zone PGF and PSP on a uniform grid within the observation domain and subsequently interpolate to the required observer points. \textcolor{black}{For a target object of size $D_x,D_y,D_z$}, we construct a sparse uniform grid of observation points \{$\mathbf{r}_{n(l,p,q)}^o={x_{l}^o, y_{p}^o, z_{q}^o}$\} (red grid in Fig.~\ref{fig:grid}) and a spatially shifted uniform grid of source points \{$\mathbf{r}_{n(l,p,q)}^s={x_{l}^s, y_{p}^s, z_{q}^s}$\} (black grid in Fig.~\ref{fig:grid}):
\begin{subequations}\label{eq:grid}
    \begin{align}
        &x_{n(l,p,q)}^s = \frac{(l-1)(D_x-\Delta_x)}{N_{gx}-1}, y_{n(l,p,q)}^s = \frac{(p-1)(D_y-\Delta_y)}{N_{gy}-1}, z_{n(l,p,q)}^s = \frac{(q-1)(D_z-\Delta_z)}{N_{gz}-1}, \label{eq:grida} \\ 
        &x_{n(l,p,q)}^o = x_{n(l,p,q)}^s + \frac{\Delta_x}{2},\ \ \ \ \  y_{n(l,p,q)}^o = y_{n(l,p,q)}^s + \frac{\Delta_y}{2},\ \ \ \ \ \ \ z_{n(l,p,q)}^o = z_{n(l,p,q)}^s + \frac{\Delta_z}{2}, \label{eq:gridb} \\
        &\Delta_x = \frac{D_x}{N_{gx}},\ \Delta_y = \frac{D_y}{N_{gy}},\ \Delta_z = \frac{D_z}{N_{gz}};\ n= l N_{gz}N_{gy} + p N_{gz} + q. \label{eq:gridc}
    \end{align}
\end{subequations}
Here, $N_{gx}$, $N_{gy}$, $N_{gz}$ are the number of the grid points in the $x$, $y$, $z$ dimensions, $N_g=N_{gx}N_{gy}N_{gz}$ is the total number of the grid points, which is of  $O(1)$,  with $n\in [1, N_g]$, and $\Delta_x$, $\Delta_y$, $\Delta_z$ are grid spacing on  in the $x$, $y$, $z$ dimensions. The source grid is shifted from the observer grid by $\Delta_x/2$, $\Delta_y/2$, $\Delta_z/2$ to make the grids non-overlapping, which results in a fast convergence of the sum in Eq.~\eqref{eq:pgf} and Eq.~\eqref{eq:lgf} when calculating PGF. 

Based on this grid, the PGF at any observation point $\mathbf{r}_m$ from any source point $\mathbf{r}_n$ in the zeroth unit cell can be calculated by interpolation: 
\begin{equation}\label{eq:PGF_interp}
    G^p_{far}(\mathbf{r}_m -\mathbf{r}_n)=\sum_{m'=1}^{N_{go}(\mathbf{r}_{m})}\sum_{n'=1}^{N_{gs}(\mathbf{r}_{n})}\omega^o(\mathbf{r}_m, \mathbf{r}_{m'}^o)G^p_{far}(\mathbf{r}_{m'}^o-\mathbf{r}_{n'}^s)\omega^s( \mathbf{r}_{n'}^s,\mathbf{r}_n).
\end{equation}
Here, $\omega^s( \mathbf{r}_{n'}^s,\mathbf{r}_n)$ represents the interpolation coefficients from the source grid points $\mathbf{r}_{n'}^s$  to original non-uniform source points $\mathbf{r}_{n}$ and $\omega^o(\mathbf{r}_{m}, \mathbf{r}_{m'}^o)$ represents the interpolation coefficients from the observer grid points $\mathbf{r}_{m'}^o$ to the original non-uniform observer points $\mathbf{r}_{m}$. The number of the source interpolation points $N_{gs}(\mathbf{r}_{n})$ and observer interpolation points $N_{go}(\mathbf{r}_{n})$ represents the number of the grid points that need to be used for the interpolation. For example, we can choose the $q$-th order Lagrange interpolation such that all grid points participate in the interpolation for any of the non-uniform points. For this choice, $N_q = (q+1)^3$ , \textcolor{black}{$N_{go}=N_{gs}=N_q$}, and the coefficients $\omega^o$ and $\omega^s$ are the same for the same points. For a more general case of a grid with $N_q$ greater than $N_{go}$ and $N_{gs}$ , the interpolation for each non-uniform point involve only a subset of the grid points, and the coefficients $\omega^o$ and $\omega^s$ can be represented as sparse matrices. The corresponding sparse matrices are generally transpose versions of each other.

We can substitute $G^p_{far}(\mathbf{r}_m - \mathbf{r}_n)$ from Eq.~\eqref{eq:PGF_interp} into Eq.~\eqref{eq:1} and re-arrange the sums to obtain an alternative representation:
\begin{equation}
    \label{eq:PSP_interp}
    u(\mathbf{r}_m)= \sum_{m'=1}^{N_{go}(\mathbf{r}_{m})}\omega^o(\mathbf{r}_m, \mathbf{r}_{m'}^o)\sum_{n'=1}^{N_g}G^p_{far}(\mathbf{r}_{m'}^o-\mathbf{r}_{n'}^s)\sum_{n=1}^{N_s(\mathbf{r}_{n'}^s)}\omega^s( \mathbf{r}_{n'}^s,\mathbf{r}_n)q(\mathbf{r}_n).
\end{equation}

Here, $N_s$ is the number of the non-uniform source points that participate in the process of interpolations with a source grid point $\mathbf{r}_{n'}^s$. These grid points can be found from the interpolation procedure of Eq.~\eqref{eq:PGF_interp}. For example, when choosing the $q$-th order Lagrange interpolation with $N_q = (q+1)^3$, $N_s=N$. For a more general choice of a greater $N_q$, $N_s$ is found from the sparse matrix representation of the interpolation coefficients $\omega^o$ and $\omega^s$ .

Based on this grid construction in Eq.~\eqref{eq:grid} and PSP representation in Eq.~\eqref{eq:PSP_interp}, we can first pre-compute the coefficients $\omega^o$ and $\omega^s$, e.g., as a sparse matrix, and pre-compute $G^p_{far}(\mathbf{r}_{m'}^o - \mathbf{r}_{n'}^s)$ table at the grid points. The computation of $\omega^o$ and $\omega^s$ involves $O(N)$ operations since the interpolations are spatially local. The computation of $G^p_{far}(\mathbf{r}_{m'}^o-\mathbf{r}_{n'}^s)$ involves $O(N_g)$ operations since the grid is uniform and $G^p_{far}$ is shift invariant in term of the differences $\mathbf{r}_{m'}^o-\mathbf{r}_{n'}^s$. The computation of $G^p_{far}(\mathbf{r}_{m'}^o-\mathbf{r}_{n'}^s)$ can be done via the rapidly convergent sums of Eq.~\eqref{eq:pgf} and Eq.~\eqref{eq:lgf} due to the choice of the source and observer grid in Eq.~\eqref{eq:grid} with a spatial separation. The far-zone PSP $u_{far}$ can, then, be evaluated via the following three steps, which are illustrated in Fig.~\ref{fig:grid}.

\textit{Step 1, Projection from the source points to the source grid:}
This step calculates effective charges at the uniform source grid. Substituting the representation of Eq.~\eqref{eq:PSP_interp} for $G^p_{far}$ into Eq.~\eqref{eq:1} for the far-zone PGF and potential components allows obtaining an expression for the effective charge at the source grid points as projections from the original non-uniform sources (black arrows in Fig.~\ref{fig:grid} represent a contribution of a single non-uniform source to its relevant source grid points):
\begin{equation} \label{eq:projection}
q_g(\mathbf{r}_{n'}^s)=\sum_{n=1}^{N_{s}(\mathbf{r}_{n'}^s)}\omega^s(\mathbf{r}_{n'}^s,\mathbf{r}_{n})q(\mathbf{r}_n),
\end{equation}
where $n'\in [1, N_s]$. 

\textit{Step 2, Evaluation of the PSP at the observer grid:} 
This step calculates the far-zone PSP at the sparse observer grid points from the sparse source effective charges via the following convolution (green arrows in Fig.~\ref{fig:grid}):
\begin{equation}\label{eq:8}
    u_{far}^g({\mathbf{r}_{m'}^o})=\sum_{n'=1}^{N_g}G^p_{far}(\mathbf{r}_{m'}^o-\mathbf{r}_{n'}^s)q_g(\mathbf{r}_{n'}^s),
\end{equation}
which can be evaluated rapidly directly since $N_g=O(1)$.

\textit{Step 3, Interpolation from the observer grid to the observer points:}
The far-zone component of the PSP at the original non-uniform points can be finally calculated by interpolation, which can be viewed as an inverse procedure of step 1 (red arrows in Fig.~\ref{fig:grid}):
\begin{equation}\label{eq:interpolation}
    u_{far}(\mathbf{r}_m)=\sum_{m'=1}^{N_{go}}\omega^o(\mathbf{r}_{m}, \mathbf{r}_{m'}^o)u_{far}^g(\mathbf{r}_{m'}^o),
\end{equation}
where $m\in[1, N]$.

\subsection{Evaluation of the near-zone PSP} \label{ssec:3_near_field}
To evaluate the near-zone PSP component, $u_{near}$, we note that for $i_d = 0$, we only consider the zeroth order unit cell, which is equivalent to a \textcolor{black}{non-periodic} problem. The evaluation of $u_{near}$ for a free-space \textcolor{black}{non-periodic} problem can be efficiently handled using various fast computational approaches, such as the fast multipole like methods \cite{GREENGARD1987325, 6230628} or interpolation-based approaches \cite{1159856, LI20108463, MENG20108430}. For $i_d > 0$, neighbor images in the near-zone evaluation are included and they require modifying the fast evaluation algorithms. Here, we use the box-adaptive integral method of \cite{NUFFT,6348120}, which is also related to the pre-corrected FFT method \cite{662670} and adaptive integral method \cite{ali_yilmaz, 7770011}. We modify this approach to efficiently account for the inclusion of the periodic images. The idea is similar to the far-zone procedure in Sec.~\ref{ssec:3_far_field} in its first three steps, but it also adds another modified correction step.

Similar to the grid construction in Sec.~\ref{ssec:3_far_field}, we construct uniform grids. However, there are two notable differences. The first difference is that the source and observer uniform grids coincide, i.e., there is a single uniform grid \{$\mathbf{r}_m={x_m, y_m, z_m}$\} =
\{$\mathbf{r}_m^o={x_m^o, y_m^o, z_m^o}$\} = \{$\mathbf{r}_n^s={x_n^s, y_n^s, z_n^s}$\} by  removing $\Delta_x$, $\Delta_y$, $\Delta_z$ in Eq.~\eqref{eq:grida} and Eq.~\eqref{eq:gridb}. The reason for using a single grid is that the near-zone PGF only involves a small number of terms in its evaluation and does not need to have a spatial separation between source and observer points. The second difference is that the number of the grid points $N_g$ is significantly larger than what is used in the far-zone evaluation. Typically, $N_g$ is comparable to the number of the non-uniform points $N$ for the optimal performance. This large $N_g$ is required to have $N_s$ (as defined in Eq.~\eqref{eq:PSP_interp}) to be of $O(1)$ to result in a reduced number of computations in the correction step 4. Based on this grid definition, the near-zone evaluation proceeds in the following four steps.

\textit{Step 1, Projection from the source points to the source grid:}
This step is similar to step 1 in Sec.~\ref{ssec:3_far_field}, where $G^p_{far}$ is replaced with $G^p_{near}$. Here, $N_g$ is typically much larger than $N_s$, and the matrices of the interpolation coefficients are sparse. The resulting computational cost is of $O(N_g)=O(N)$.

\textit{Step 2, Evaluation of the PSP at the observer grid:}
This step is similar to step 2 in Sec.~\ref{ssec:3_far_field}, where $G^p_{far}$ is replaced with $G^p_{near}$. However, because $N_g$ is large, the convolution sum of Eq.~\eqref{eq:8} is computed via FFT, which is allowed because the grids are uniform. Since $G^p_{near}$ is not periodic, using FFT requires making the equivalent matrix circular, which involves doubling its size in each periodicity dimension. For performing FFTs in the computation stage, $G_{near}^p$ is tabulated in the pre-processing stage. The tabulation computational cost is of $O(N_g)=O(N)$ and the FFT computational cost is of $O(N \log N)$

\textit{Step 3, Interpolation from the observer grid to the observer points:}
The interpolation step for $u_{near}$ is identical to step 3 in Sec.~\ref{ssec:3_far_field} for $u_{far}$ and it can be viewed as an inverse of step 1 using the same sparse matrix for the interpolation coefficients. The resulting computational cost is of $O(N_g)=O(N)$.

\textit{Step 4: Error correction:}
The ability to perform the projection/interpolation in steps 1 and 3 rely on the assumption of slow variations of the Green's function. For the far-zone evaluation, $G^p_{far}$ indeed varies slowly in the domain of the zeroth unit cell because it excludes the direct and near-zone image interactions between the sources and observers. The near-zone component with $G^p_{near}$, however, involves direct and near-image interactions, which result in rapid spatial variations, e.g., the rate of such spatial variations is unbounded when the source and observer coincide. Therefore, the interpolations for the nearby source-observer pair contributions are highly inaccurate, and they need to be corrected. To eliminate these errors, we define an error correction range $D^{ER}$ for each box \textcolor{black}{(defined by $\Delta_x$, $\Delta_y$ and $\Delta_z$ in Eq.~\eqref{eq:gridc})}, determined by the uniform grid. This range is of the same order as $max\{\Delta_x,\Delta_y,\Delta_z\}$. For all observer points in a box, the corrections are performed in the error-correction region, $\Omega^{ER}(\mathbf{r}_m)$, for all the non-uniform sources in the same box and a certain number of surrounding boxes that are withing the $D^{ER}$ distance from the box. The correction involves subtracting the contributions due to steps 1-3 and adding the exact point-to-point superposition sums via appropriate contributions in Eq.~\eqref{eq:1}. The error correction procedure can be expressed as

\begin{subequations} \label{eq:er}
    \begin{align}
        &u_{near}(\mathbf{r}_m) = u_{near}(\mathbf{r}_m) + \sum_{n}^{\mathbf{r}_n \in \Omega^{ER}(\mathbf{r}_m)}[G^p_{near}(\mathbf{r}_m-\mathbf{r}_n) - G_{grid}(\mathbf{r}_m- \mathbf{r}_n)]q(\mathbf{r}_n), \label{eq:er_1}\\
        &G_{grid}(\mathbf{r}_m, \mathbf{r}_n) = \sum_{m'}^{\omega^o(\mathbf{r}_m,\mathbf{r}_{m'})\neq 0}\:\sum_{n'}^{\omega^s(\mathbf{r}_{n'},\mathbf{r}_n)\neq0}\omega^o_{nufft}(\mathbf{r}_m,\mathbf{r}_{m'})G^p_{near}(\mathbf{r}_{m'}- \mathbf{r}_{n'})\omega^s_{nufft}(\mathbf{r}_{n'},\mathbf{r}_n), \label{eq:er_2}
    \end{align}
\end{subequations}
where the error is corrected in the region $\Omega^{ER}(\mathbf{r}_m)$, i.e., for any point $\mathbf{r}$ within this region, $|\mathbf{r}-\mathbf{r}_m| \leq D^{ER}$. Similar to Eq.~\eqref{eq:projection} and Eq.~\eqref{eq:interpolation}, the $\omega^s_{nufft}$ and $\omega^o_{nufft}$ correspond to the interpolation matrices containing the weights for the above step 1 and 2. The function $G_{grid}$ represents the near-zone Green's function between the source and observer points obtained via the grid interactions with interpolation, which needs to be subtracted in the error correction process.
While the procedure in Eq.~\eqref{eq:er} can be efficient, it still can significantly increase the computational cost due to the need to account for all the near-zone images via $G^p_{near}$. For example, in the 3D periodicity case, when including 2 surrounding boxes per each dimension \textcolor{black}{($i_d=2$), there are 125 ($(2\times i_d+1)^3$) images in total}, which leads to a high additional computational cost. 

\begin{figure}[h] 
\includegraphics[width=10cm]{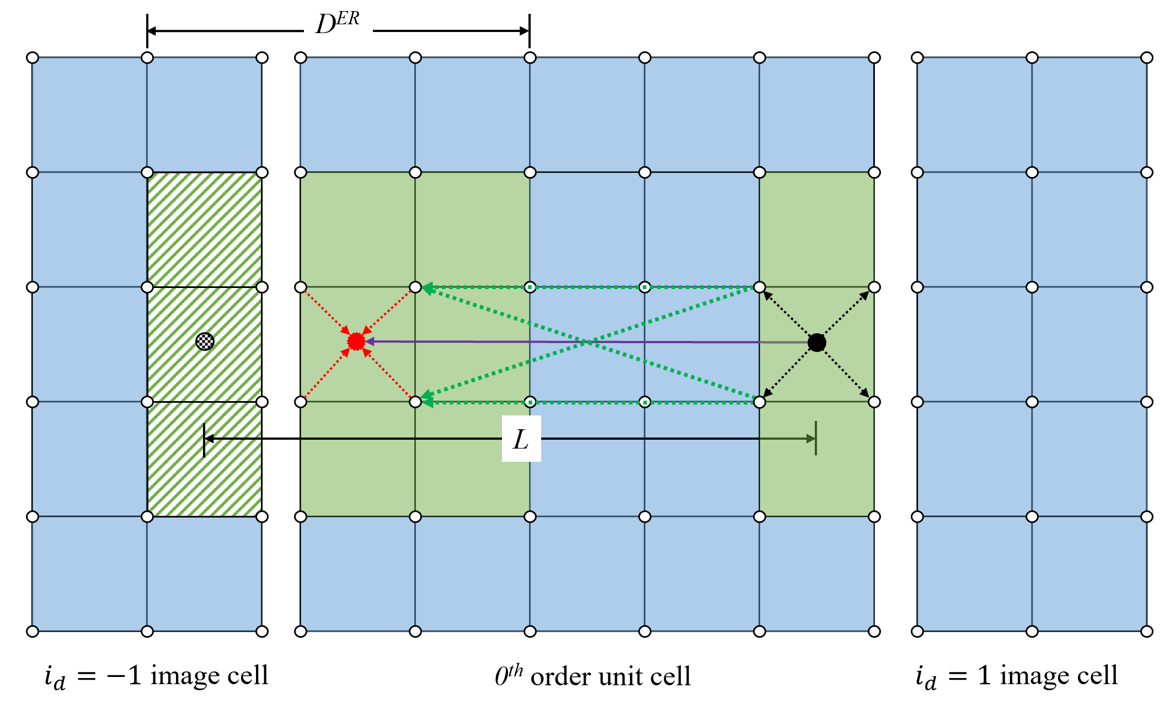}
\centering
\caption{Illustration for the error correction step 4 for the near-zone PSP component. The black dot is a source point and red dot is an observer point. The green region is the error-correction region $\Omega^{ER}$ . The black circles are grid points. For the observers in $\Omega^{ER}$, PSPs are calculated through direct calculations, this is done by subtracting the grid interaction inside the green region (dotted arrows) and adding the direct interactions (dashed purple arrow). The shaded green region and grey dot is the image of the right boundary in the $i_d=-1$ image cell, which is close to the left boundary of the unit cell within $\Omega^{ER} $.} \label{fig:nufft}
\end{figure}

This computational cost can be much reduced by including only a single extra periodic image per periodicity dimension. This is possible due to the fact that $\{\Delta_x,\Delta_y,\Delta_z\}$ are typically much smaller than $\{D_x,D_y,D_z\}$, so that most of the images correspond to a large spatial separation. For example, in Fig.~\ref{fig:nufft} showing an example for a 1D $x$-dimension periodicity case, for the sources in the internal boxes of the zeroth (middle) unit cell, only the $G_0(\mathbf{r}_m-\mathbf{r}_n)$ component of the $G^p_{near}(\mathbf{r}_m-\mathbf{r}_n)$ can be used for the near-zone correction. On the other hand, for the left box, which is next to the left edge of the zeroth unit cell (the left green box with red dot and arrows), in addition to the $G_0(\mathbf{r}_m-\mathbf{r}_n)$ component, the component of $G_0(\mathbf{r}_m-\mathbf{r}_n-L_x\mathbf{\hat{x}})$ of $G^p_{near}(\mathbf{r}_m-\mathbf{r}_n)$ needs to be accounted for. The latter component corresponds to the sources from the right box (right green box with a black dot and arrows) contributing to the left box observers via the $i_x=1$ image. These sources are equivalent to the sources in the right box of the left image cell (the dashed green box with a grey dot) contributing to the left box observers via the image Green's function $G_0(\mathbf{r}_m-\mathbf{r}_n)$, i.e., they are equivalent to the geometrically close boxes. This procedure is generalized to 2D and 3D periodicities by including corresponding required images. The computational cost of this procedure is only slightly higher than that for a \textcolor{black}{\textcolor{black}{non-periodic}} case, in which only the $G_0(\mathbf{r}_m-\mathbf{r}_n)$ contribution is required. The cost increase is associated with the need to include additional images at the sources associated with the boxes located at the corners, edges, and sides of the computational domain. There is a small number of such boxes as compared to the total number of near-zone images, so that the computational cost increase is insignificant. The overall computational cost of this step is of $O(N)$.

\section{Numerical results} \label{sec:4}
We implemented FFT-PIM both in a Central Processing Unit (CPU) and Graphics Processing Unit (GPU) based code. This section presents numerical test results demonstrating the accuracy and performance of the FFT-PIM. Sec.~\ref{ssec:4.1} presents error analysis, comparing the results from the FFT-PIM with results obtained via the direct sum of Eq.~\eqref{eq:1} with Eq.~\eqref{eq:lgf_1} for different problem sizes, uniform grid sizes, and interpolation orders. Sec.~\ref{ssec:4.2} examines the computational speed of the FFT-PIM for different parameters, such as problem size, uniform grid size, and the periodicity dimension.

\begin{figure}[htbp] \label{fig:potential}
    \centering
    \subfigure[]{
    \includegraphics[width=0.47\linewidth]{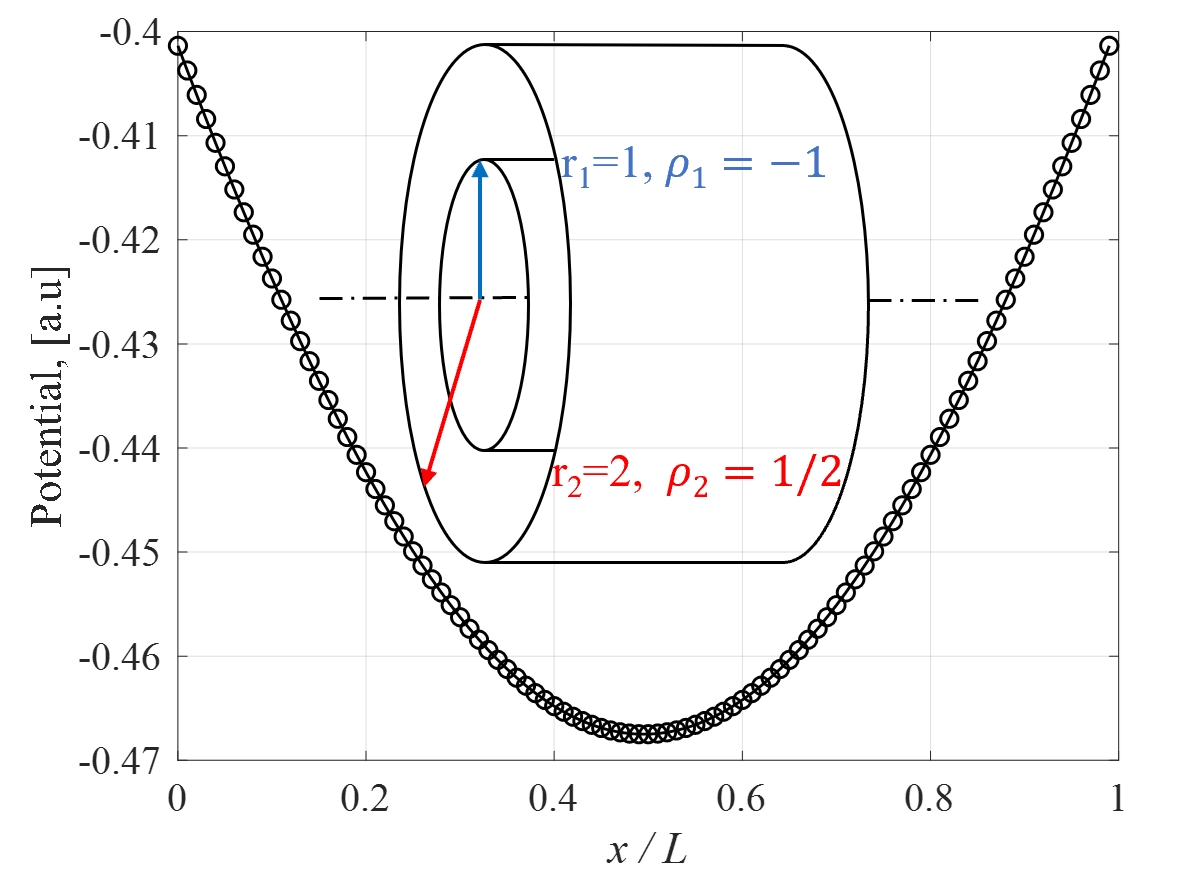}\label{fig:potentiala}
    }
    \quad
    \subfigure[]{
    \includegraphics[width=0.47\linewidth]{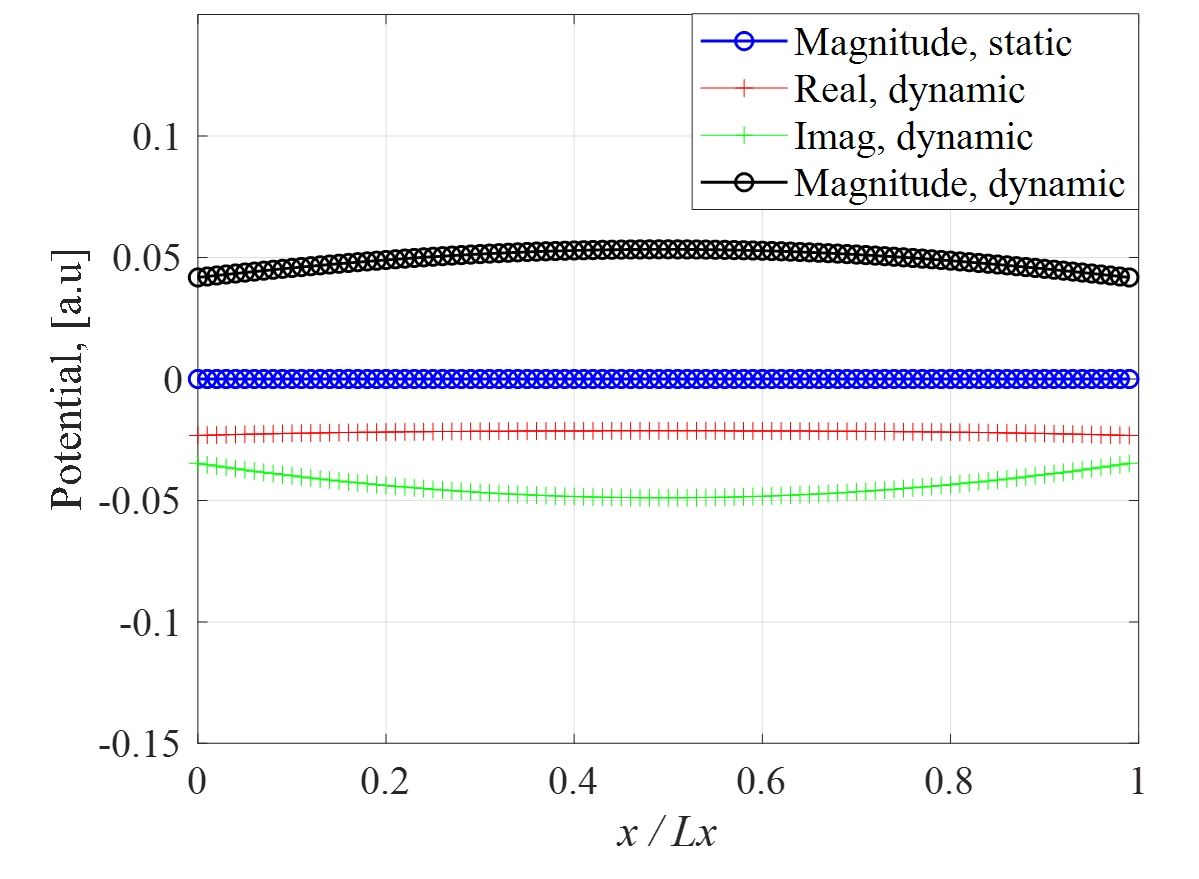}\label{fig:potentialb}
    }
    \caption{Potential on the $x$ axis of a coaxial structure . (a) Illustration of coaxial structure unit cell ($L=1$) and potential of the \textcolor{black}{non-periodic} unit cell. Inner radius $r_1=1$ with negative line charge density $\rho_1=-1$, outer radius $r_2=2$ with positive line charge density $\rho_2=1/2$. (b) PSP with a 1D periodicity along the $x$ axis with $L_x$ = 1. The blue line is the magnitude of PSP for the \textcolor{black}{NPSP} case with $k_0=k_{x0}=0$. The red, green and black lines are real part, imaginary part, and magnitude of PSP for the dynamic case with $k_0=1$, $k_{x0}=1-j$. } \label{fig:potential}
\end{figure}

We start by showing results for the calculated potential due to a source distribution (Fig.~\ref{fig:potential}). The structure is 1D periodic in the $x$-direction with the periodicity of $L_x = 1$ and the unit cell comprises a section of coaxial cable along the $x$-axis that extends through the entire period. The inner and outer radius is $r_1=1$ and $r_2=2$, respectively. Negative and positive sources (e.g., charges) are distributed uniformly on the inner surface with charge density $\rho_1=-1$ and $\rho_2 = 1/2$, respectively, resulting in a neutral unit cell. We first calculate the potential on the central $x$ axis for a \textcolor{black}{\textcolor{black}{non-periodic}} case, as shown in Fig.~\ref{fig:potentiala} and we observe significant edge effects in terms of the potential spatial dependence. We then use the FFT-PIM with a 1D periodicity, making it into an infinitely long coaxial source distribution. For the \textcolor{black}{NPSP} case, we find that the potential becomes near zero everywhere (blue curve in Fig.~\ref{fig:potentialb}), the edge effect is eliminated, which agrees with the analytical result. We also calculate the potential for the dynamic case with the wavenumber $k_0=1$ and a complex phase shift $k_{x0}=1-j$, in which case the potential becomes complex with non-zero magnitude and corresponding real/imaginary parts (black, red, and green curves in Fig.~\ref{fig:potentialb}).

We, then, proceed with showing results for the error analysis and computational performance.

\subsection{Error analysis} \label{ssec:4.1}
To test the error of the FFT-PIM, we use an infinite long 1D bar along the $x$ axis, and the target object is a cube of the size $D_x,D_y,D_z=50$ with 1D $x$-direction periodicity with $L_x=D_x=50$. We assume a \textcolor{black}{NPSP} case, i.e., $k=k_{x0}=0$. \textcolor{black}{We mesh the cube with a tetrahedral mesh and the total number of the mesh points (vertices) is $N=7189$. In order to calculate the error, we randomly set 7000 mesh points with finite values as source points while the rest 189 points are set to $0$, and make sure that $\sum_{l=1}^N q_l = 0$. The 189 mesh points set to $0$ are considered as observer points, and the locations of these observer points are chosen such that they are spatially shifted as compared to the source points such that the PGF in Eq.~\eqref{eq:lgf_1} converges for all source-observer pairs when the direct superposition sum is used.} We compare the results for relative errors between the FFT-PIM and the direct computation via Eq.~\eqref{eq:lgf_1}, which is equivalent to Eq.~\eqref{eq:pgf_1} as shown in Sec.~\ref{sec:2}. We tested the relative errors at the observer points with different interpolation orders, uniform grid size, and number of subtracted near-zone images. The results are shown in Fig.~\ref{fig:error}. We find that with cubic interpolation order and a uniform grid of size $N_{gx},N_{gy},N_{gz}=10$, by subtracting the zeroth order and only the nearest neighbors ($i_d = 1$), we can achieve an error at the level of $10^{-3}$, which is sufficient for many practical cases. The relative error is further reduced by increasing the number of grid points, order, and the near-zone images.

\begin{figure}[htbp] \label{fig:error}
    \centering
    \subfigure[]{
    \includegraphics[width=0.47\linewidth]{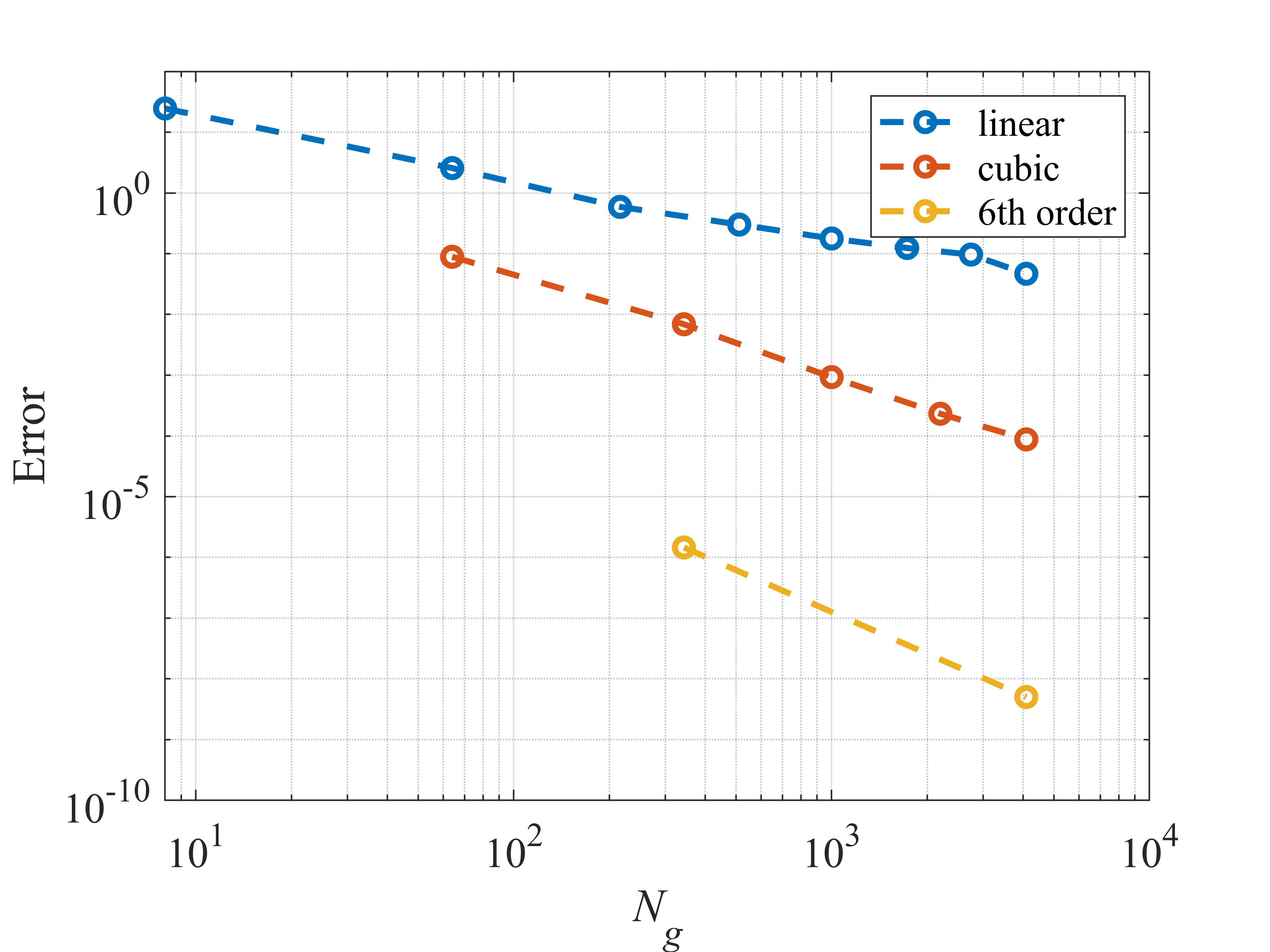} \label{fig:error_1}
    }
    \quad
    \subfigure[]{
    \includegraphics[width=0.47\textwidth]{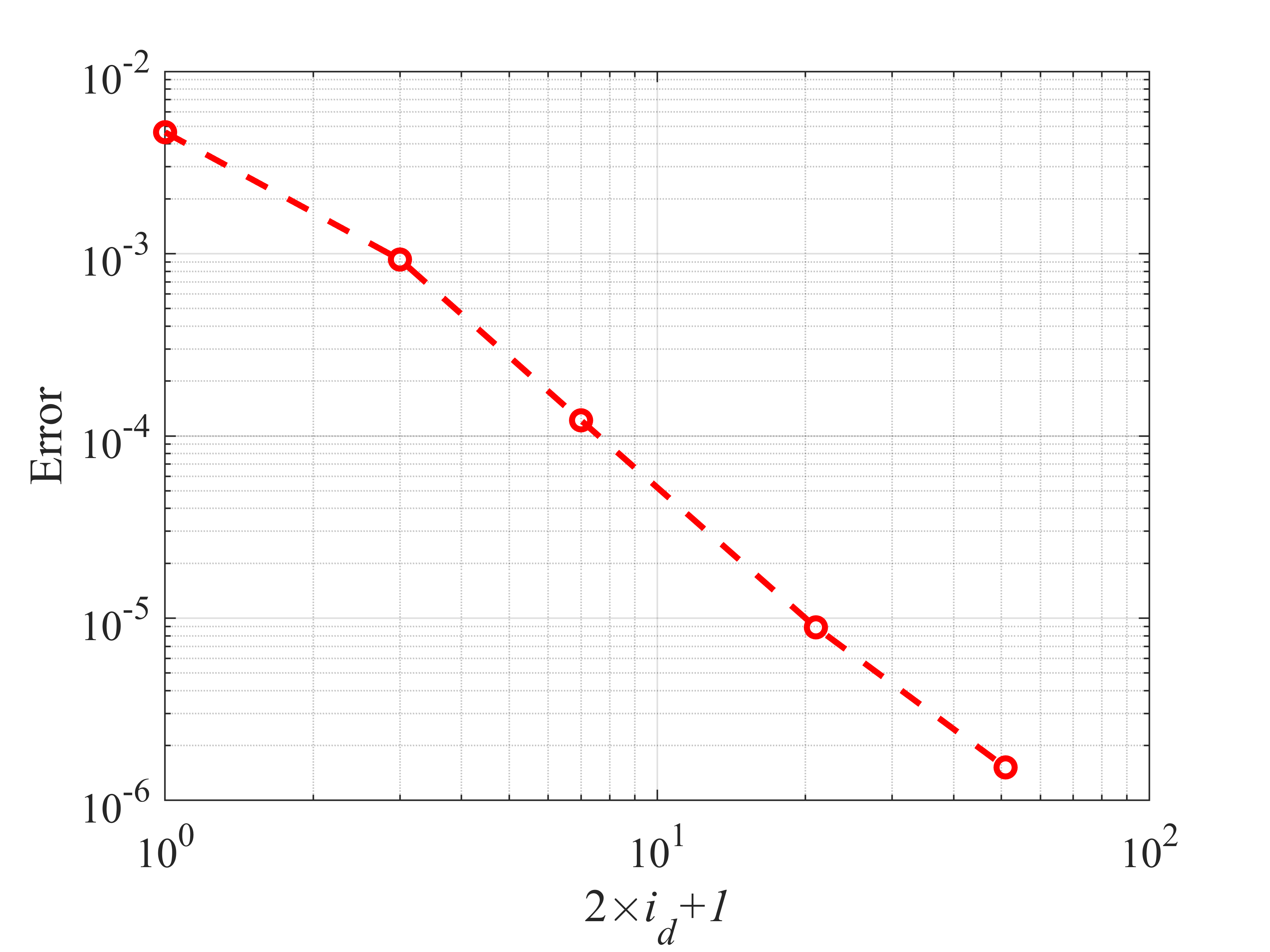} \label{fig:error_2}
    }
    \caption{Relative error of the far-zone PSP for a \textcolor{black}{NPSP} case with a 1D periodicity along the $x$-axis. (a) Relative error for 1st, 3rd, and 6th order interpolation with varying uniform grid sizes with $i_d=1$. (b) Relative error of cubic interpolation with 1000 sparse grid points and different numbers of subtracted near-zone unit cells.} \label{fig:error}
\end{figure}

\subsection{Computational performance analysis} \label{ssec:4.2}

\begin{figure}[htbp] 
\includegraphics[width=8cm]{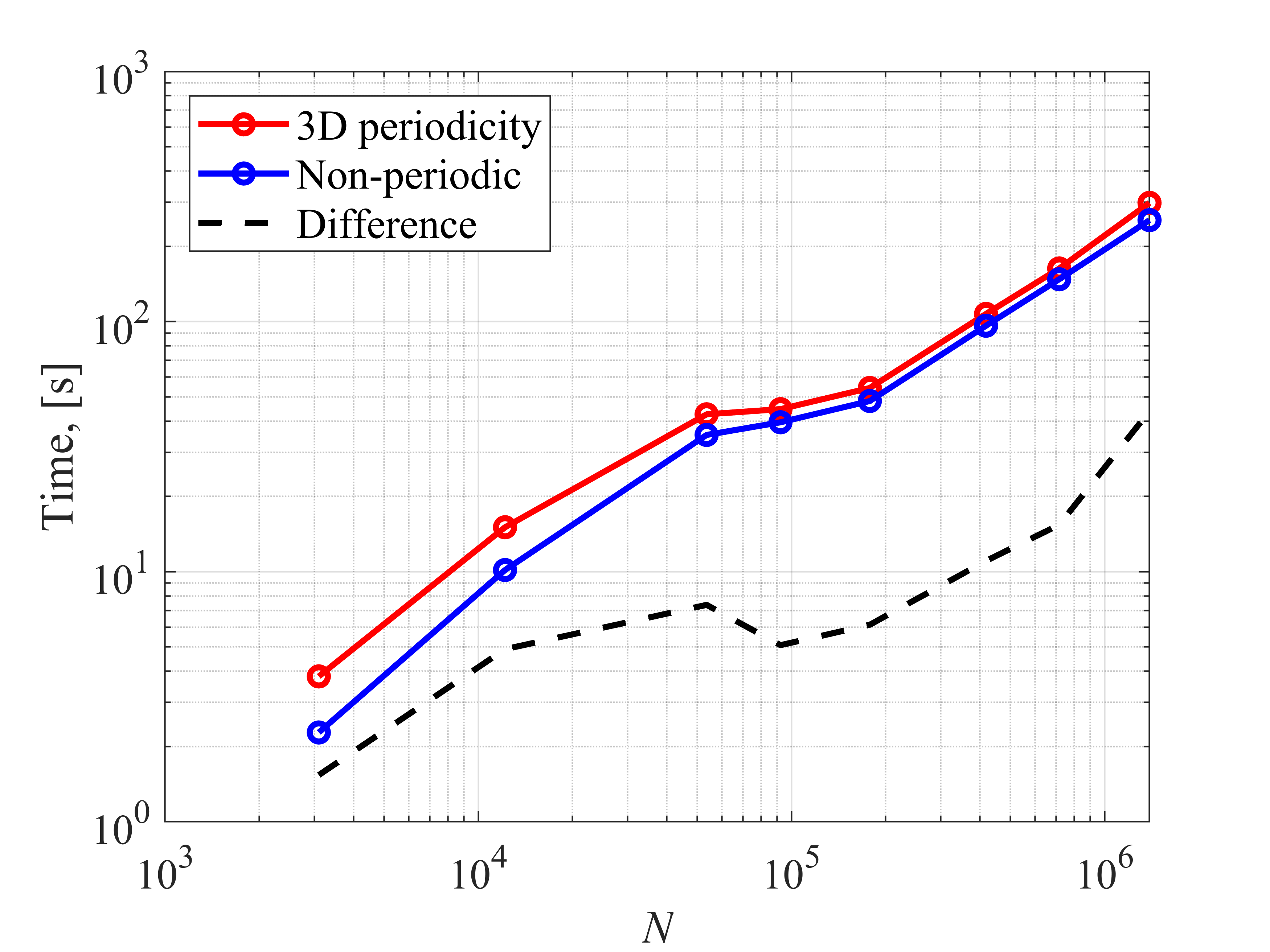}
\centering
\caption{\textcolor{black}{Preprocessing time for the non-periodic case, NPSP case with 3D periodicity, and for their difference.}} \label{fig:preprocess}
\end{figure}

\begin{figure}[htbp]
    \centering
    \subfigure[]{
    \includegraphics[width=0.47\linewidth]{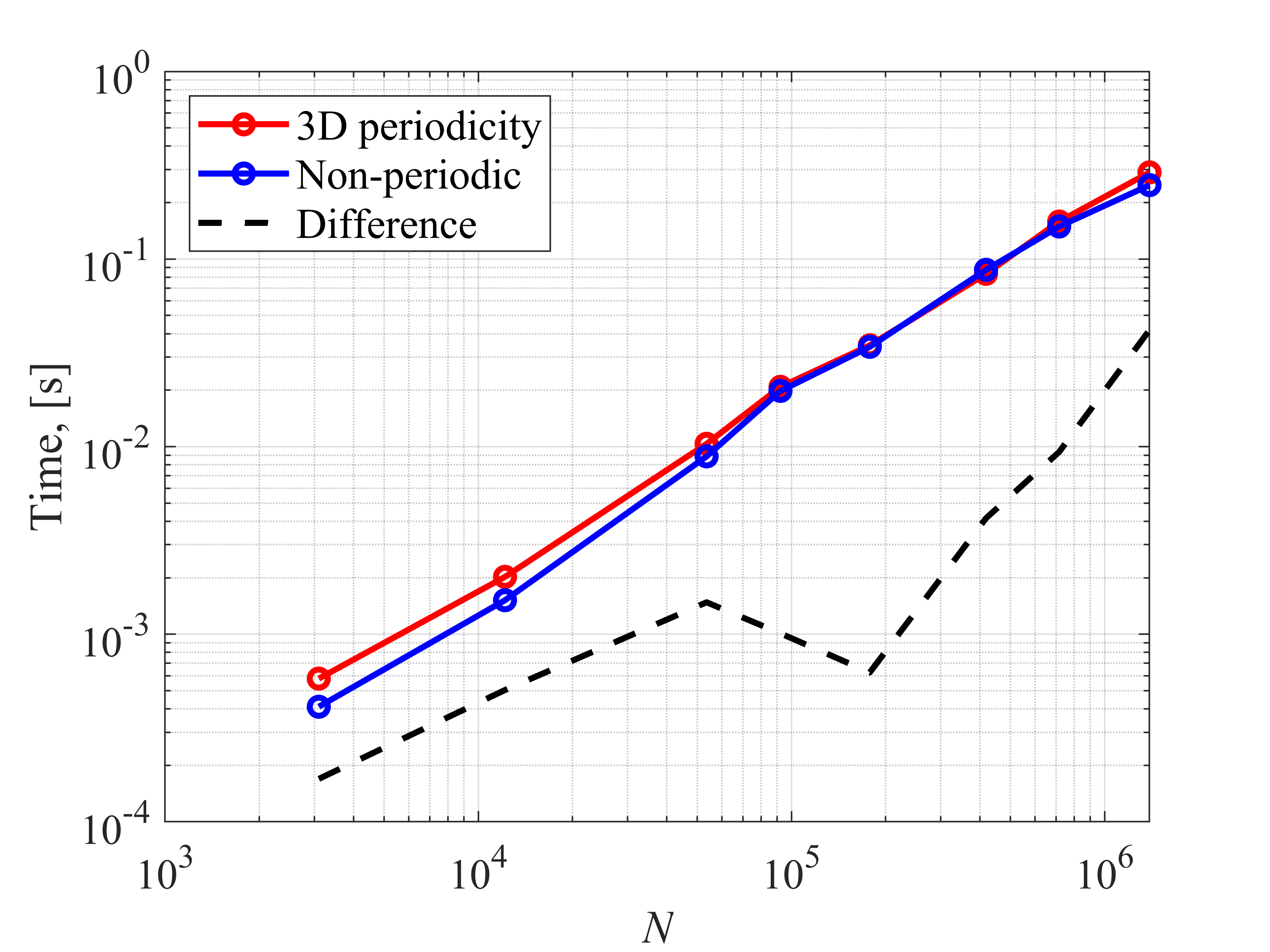}
    }
    \quad
    \subfigure[]{
    \includegraphics[width=0.47\textwidth]{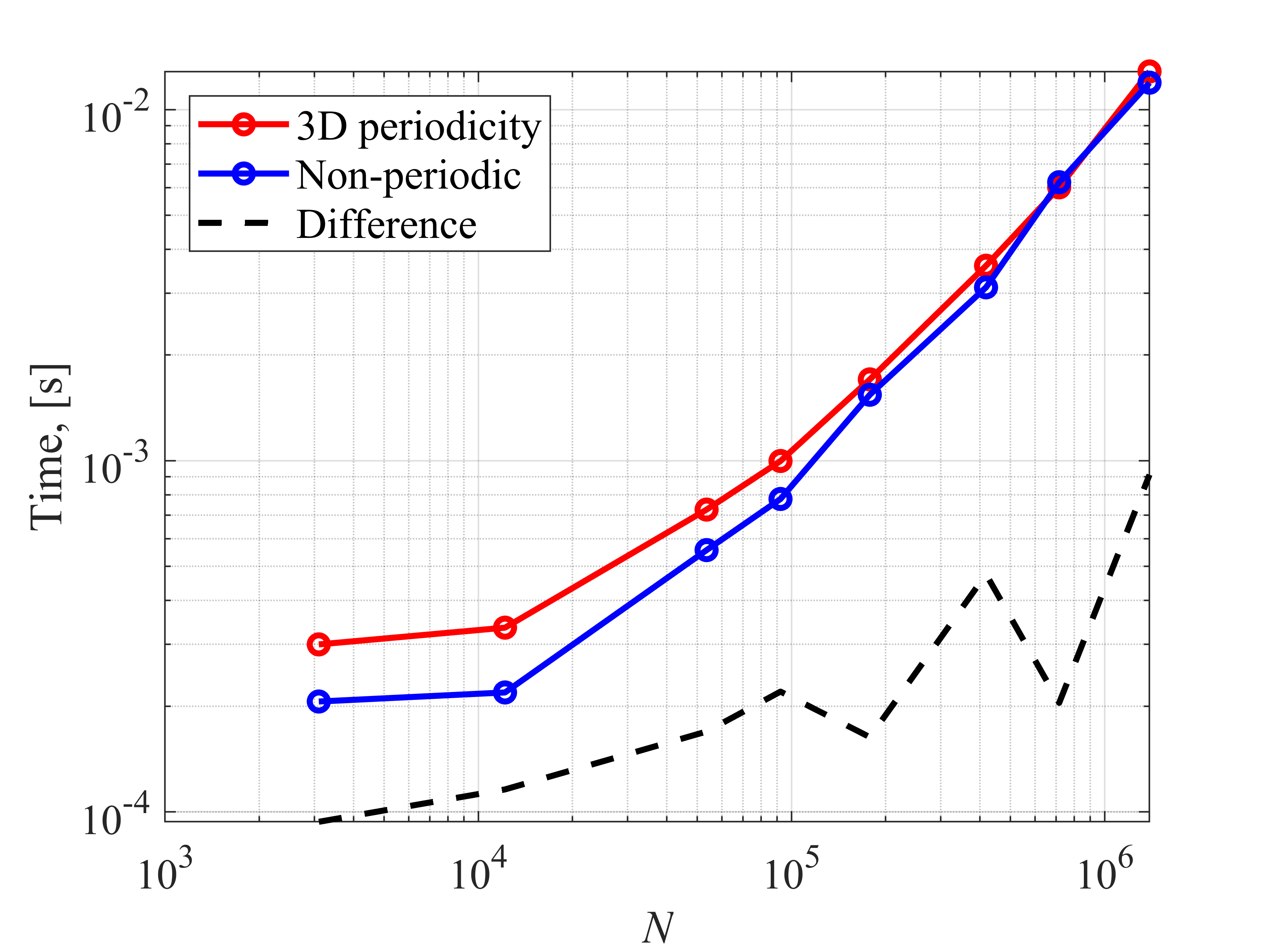}
    }
    \caption{\textcolor{black}{Execution time versus $N$ of evaluating the non-periodic potential and \textcolor{black}{NPSP} PSP with 3D periodicity on (a) eight-core CPU and (b) on GPU.}} \label{fig:performance}
\end{figure}

\begin{table}[htbp]
  \centering
  \begin{tabular}{ccc}
    \hline
    $N$ & Near-zone grid size   \\
    \hline
    3096 & $17^3$\\
    12175 & $25^3$ \\
    53601 & $41^3$ \\
    92233 & $49^3$ \\
    177973 & $59^3$ \\
    418308 & $77^3$ \\
    1391742 & $111^3$
\end{tabular}
\caption{\textcolor{black}{Grid size for near-zone PSP component evaluation. The near-zone grid is much larger than the far-zone grid when the problem size $N$ is large. Second-order projection/interpolation for near-zone evaluation is used.}} \label{tab:near-zone}
\end{table}

\begin{figure}[htbp] 
\includegraphics[width=8cm]{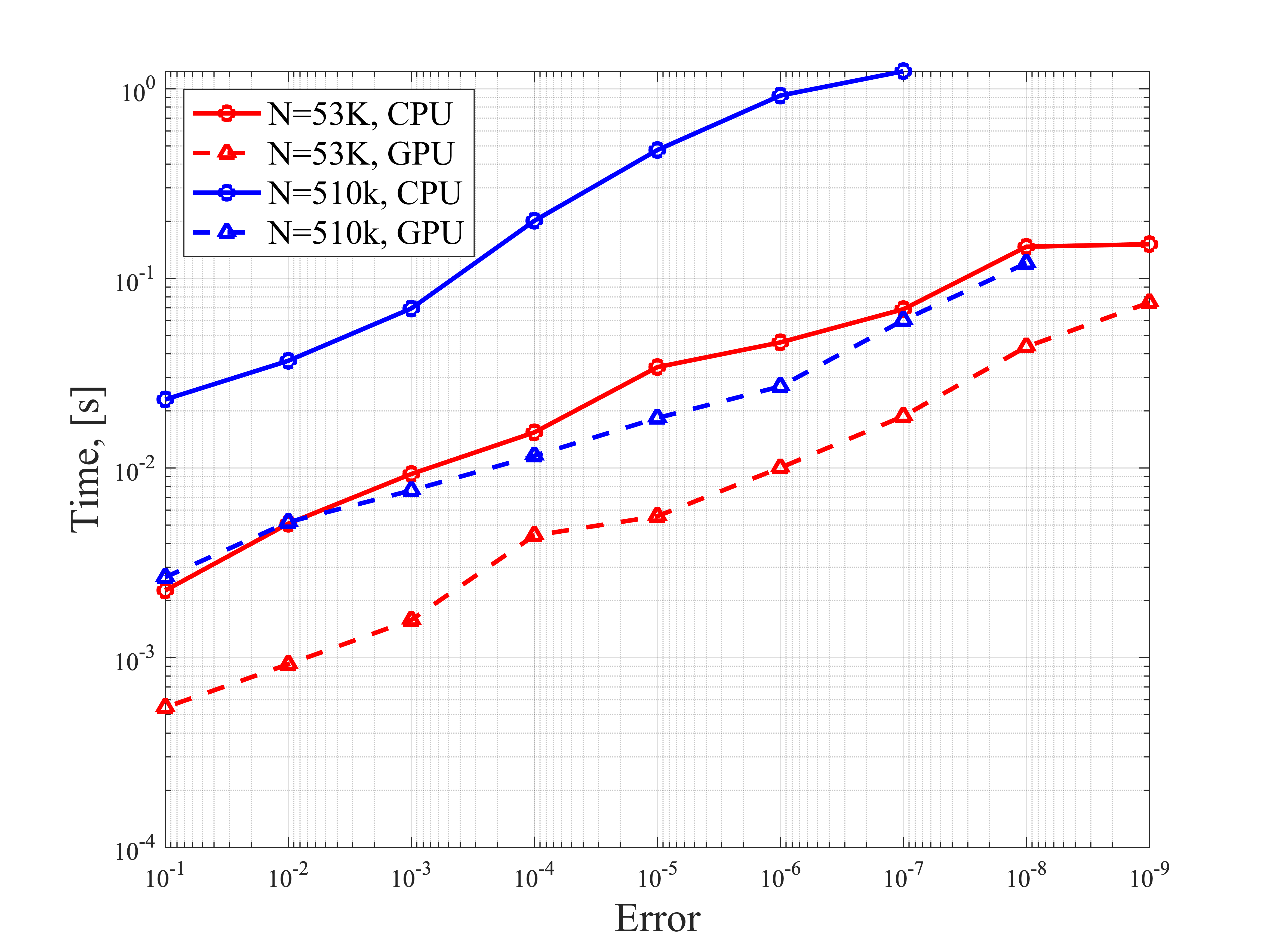}
\centering
\caption{\textcolor{black}{Execution time versus the relative error for the near-zone PSP component on CPU and GPU for a \textcolor{black}{3D NPSP} case with $N=53K$ and $N=510K$.}} \label{fig:error_perf}
\end{figure}

To benchmark the performance of the FFT-PIM, \textcolor{black}{we conducted a set of tests on Purdue
Anvil cluster at Rosen Center for Advanced Computing (RCAC) in Purdue University. CPU tests were run on eight cores of AMD Epyc "Milan" processors and GPU tests were run on Nvidia A100.} The results are shown for a 3D  \textcolor{black}{NPSP} problem. The target object is a cube with dimensions $D_x,D_y,D_z=100$ and $L_{x,y,z}=D_{x,y,z}+1$. This problem has $L_{x,y,z}\approx D_{x,y,z}$ and it represents the maximal computational complexity as compared to 1D and 2D periodicity cases. The configuration of the far-zone calculation is set as cubic interpolation with uniform grid size of $N_{gx},N_{gy},N_{gz}=10$, and $i_d = 1$. The near-zone is handled using the procedures from the FastMag micromagnetic simulator \cite{Fastmag, 6348120} for \textcolor{black}{non-periodic} and periodic cases. When periodicity presents, it is modified to allow for an efficient handling of periodicities, as outlined in Sec.~\ref{ssec:3_near_field}.

Fig.~\ref{fig:preprocess} shows the computational time on CPU for pre-processing needed for computing the interpolation coefficients and PGF tables. We compare the CPU time with the case of the same source - observer distributions but without any periodicity, which demonstrates the overhead due to the periodicity. The predominant additional time complexity arises from substituting the original free-space Green's function, $G_0$, with a more complicated $G_{near}$. Specifically in this 3D case with $i_d=1$, $G_{near}$ is $(2\times1+1)^3$ times slower than $G_0$. However, the overall discrepancy is marginal, resulting in only a 15\% overhead. This small impact is due to the fact that tabulating the Green's function constitutes only a small part of the pre-processing time, and the time complexities of other segments do not increase as the periodicity is introduced. We also evaluated 1D and 2D periodic cases, and obtained similar performance with a reduced overhead.

Fig.~\ref{fig:performance} shows the execution time for a single evaluation of PSP on CPU and GPU. When compared to a \textcolor{black}{non-periodic} problem, the major execution time increase is due to the evaluation of the far-zone PSP component, which can be performed separately from the near-zone evaluation. The presence or absence of periodicity has minimal impact on the near-zone evaluation since all Green's function values are tabulated during the pre-processing step. Consequently, in the eight-core CPU benchmark, only a minor difference ($< 15\%$) is observed between configurations with and without periodicity for moderately large case where $N > 10^4$. Furthermore, in the GPU benchmark, NVIDIA CUDA multi-stream concurrency enables simultaneous implementation of both near-zone and far-zone evaluations. This approach effectively eliminates the overhead induced by the far-zone evaluation, which is already marginal in eight-core CPU execution, resulting in \textcolor{black}{a small differences ($< 5\%$)} between cases with and without periodicity, including the 3D periodicity cases. \textcolor{black}{The execution performance on GPU is generally 20-30X faster than that on eight-core CPU for the tested sizes, which is equivalent to around $\sim$200X speedup between GPU and single CPU core.} The 1D and 2D periodicity cases exhibit a similar performance, indicating a highly efficient implementation for all 1D, 2D, 3D periodicity types. \textcolor{black}{To further present differences between the near-zone and far-zone component, we show the grid size of the near-zone PSP component in Table~\ref{tab:near-zone}. For most practical problem sizes, the near-zone grid in Table~\ref{tab:near-zone} is much denser than the far-zone grid.}

\textcolor{black}{From Fig.~\ref{fig:performance} we observe that the near-zone component PSP dominates the computational load for the 1D, 2D, and 3D periodicity cases. Similarly, Figure~\ref{fig:error} reveals that the error associated with the far-zone PSP component can be reduced with minimal adjustments to the parameters and limited overhead. This observation suggests that the near-zone PSP component is a critical factor in determining the balance between overall performance and error rates. To further explore this relationship, we conducted a series of benchmarks with a cube meshed in its volume, containing $N=53,601$ points for a coarser mesh and $N=510,643$ for a finer mesh. The results are shown in Fig.~\ref{fig:error_perf}. The performance shown here demonstrate a good performance compared to available solutions \cite{r1_2}.}

\begin{figure}[htbp]
    \centering
    \subfigure[]{
    \includegraphics[width=0.47\linewidth]{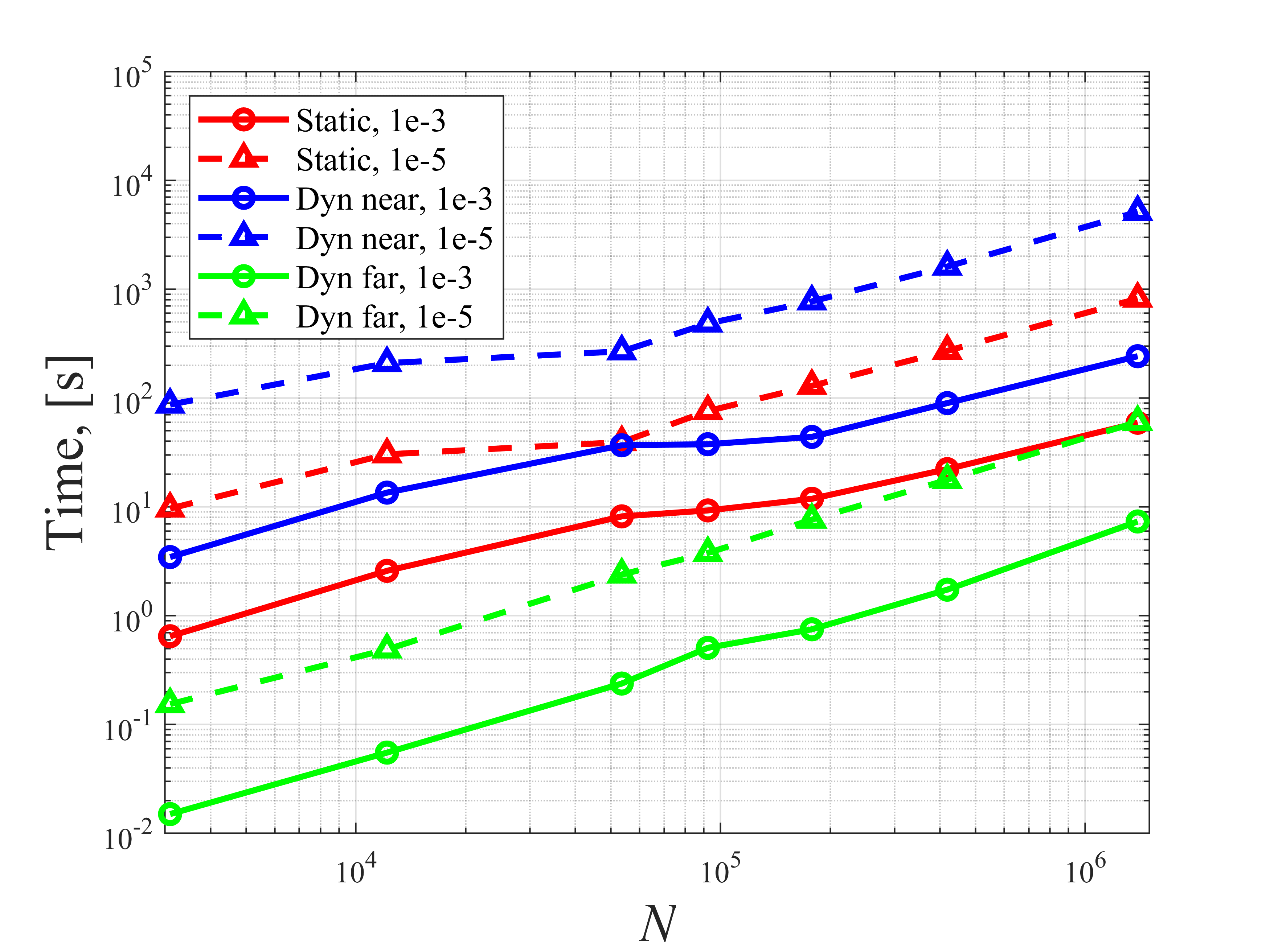}
    }
    \quad
    \subfigure[]{
    \includegraphics[width=0.47\textwidth]{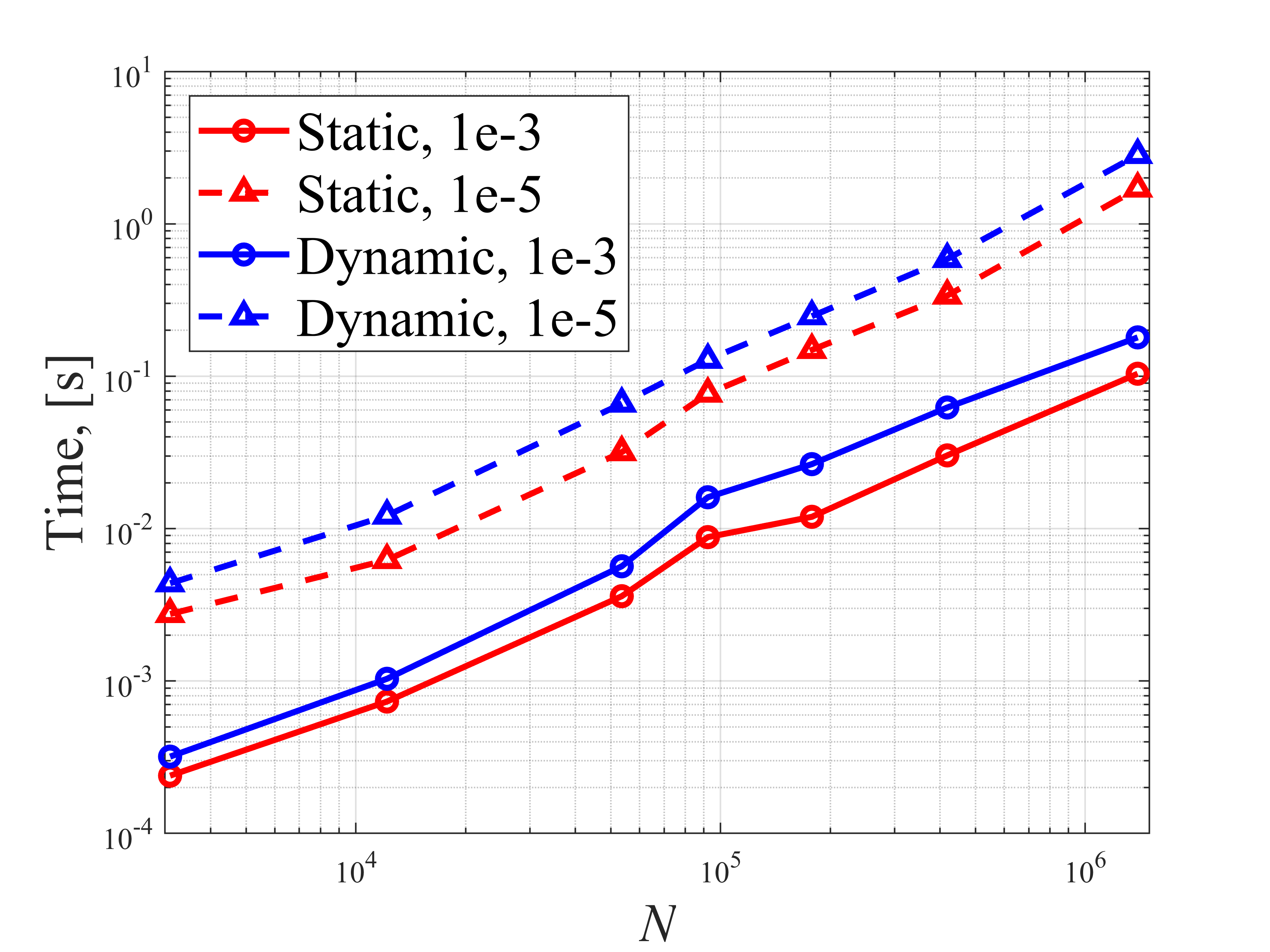}
    }
    \caption{\textcolor{black}{(a)Preprocessing time and (b) execution time versus $N$ of evaluating the PSP for a \textcolor{black}{3D NPSP} and high-frequency dynamic cases for the relative errors of $1e-3$ and $1e-5$.}}
    \label{fig:complex_real}
\end{figure}

\textcolor{black}{Figure~\ref{fig:complex_real} compares the preprocessing and computational time of the \textcolor{black}{NPSP} and dynamic cases for a 3D periodicity of $L_x =101$, $L_y =101$, $L_z =101$ for relative errors of $10^{-3}$ and $10^{-5}$. The \textcolor{black}{NPSP} case parameters are as those in Fig.~\ref{fig:performance} and Table.~\ref{tab:near-zone} for error level $10^{-3}$, and we double the near-zone grid size for error level $10^{-5}$. The dynamic cases have phase shifts of $k_{x0}=0.01-0.01j$, $k_{y0}=0.01-0.01j$, $k_{z0}=0.01-0.01j$. For the dynamic case, the wavenumber $k_0$ was chosen such that the average distance between the sources/observers was $\lambda/20$. This choice is typical when considering high-frequency dynamic problems and it allows resolving the wave process spatial variations. For the largest considered $N$, this choice results in the domain size of $L_x=L_y=L_z=5\lambda$. We find that the computational time for the low- and high-frequency cases are close to each other and they scale similarly as that for the non-periodic static case. The computational time of the dynamic cases is around 2 times greater, which is related to performing complex-valued versus real-valued operations. }

\section{Summary and discussion} \label{sec:5}
We introduced an efficient and flexible FFT-PIM for computing PSP. FFT-PIM can be used for a wide variety of problem types with the same numerical implementation. It allows computing PSP for a non-uniform source distribution, works for arbitrary 1D, 2D, and 3D periodicities, can operate with or without a phase shift between the periodic boundaries, and is applicable to dynamic and static problems. A requirement is imposed for using FFT-PIM for \textcolor{black}{NPSP} problems, which is that the source must be neutral within the unit cell. Such a requirement is natural to many practical problems, such as problems dealing with electric dipoles, magnetization, and molecular structures.

FFT-PIM is based on a superposition sum between the source distribution and PGF. PGF, which originally is defined as an infinite sum that does not have a rapid convergence, can be represented as an exponentially convergent sum in terms of spectral series expansions, which are given for the dynamic case and static case with and without a periodic phase shift. The expansions converge under the assumption of a sufficient separation between the source and observer points. To lead to a fast superposition sum computation, PGF and the corresponding PSP are represented in terms of the near- and far-zone components. The near-zone component includes a finite, typically small, number of periodic images, whereas the far-zone component includes all the rest of an infinite number of images. The far-zone component is recognized as a slowly varying function of the spatial coordinates, which allows computing the far-zone PSP component by first computing it on a sparse uniform grid, i.e., at a small number of points, and then interpolating it to all the required non-uniform observation points by local interpolation. This process involves defining shifted source and observer grids to allow for a rapidly convergence of the spectral sums for PGF. This process requires a small number of operations and low computational cost, while allowing for a rapid converge. The near-zone PSP component can be evaluated by any fast method, but it required accounting for the additional near-zone images, which may substantially increase the computational cost and memory consumption. An approach based on the box-adaptive adaptive integral method is presented that allows evaluating the superposition sum for the far-zone component rapidly based on including essential images only at the required locations near the periodic boundaries. The result is an approach that has an overall computational cost of $O(N\log N)$ and memory consumption of $O(N)$. The presented results demonstrate the high convergence, accuracy, and computational performance of the FFT-PIM.

We note that FFT-PIM can be regarded as an extension of the fast periodic interpolation method \cite{5582244} and the non-periodic box-adaptive integral method \cite{6348120} or pre-corrected FFT method \cite{662670}. With respect to the fast periodic interpolation method \cite{5582244}, the extensions are in enabling both dynamic and static cases with any phase shifts or no phase shift as well as in enabling the complete fields computation, including fast near-zone computations via an FFT-based method without a need of extending the computational domain. With respect to the non-periodic box-adaptive integral method \cite{6348120} or pre-corrected FFT method \cite{662670}, FFT-PIM allows considering problems with 1D, 2D, and 3D periodicities with almost the same computational cost as the non-periodic approaches. \textcolor{black}{We also note that the approach of separating the computations into the near- and far-zone components is related to other fast evaluation techniques \cite{662670, r1_1, r1_2, r1_3, r1_4, 5582244, 6348120} for periodic and non-periodic domains. Contributions of this work are in the efficient evaluation of the far-zone component via the properly chosen uniform grids and Green's function evaluation as well as in the adaptation of the near-zone component computation via an FFT-based method.}

While FFT-PIM is a powerful method, it also has some limitations. In particular, for high-frequency problems with electrically large periodicities, the computation of PGF may become slow. This time may be reduced by using alternative methods for computing PGF, e.g., see Refs.~\cite{6059491} and \cite{4907048}. Additionally, using FFTs mean that the computations are done on the entire FFT grid, which may be inefficient for problems of curved linear or surface domains or in cases of highly non-uniform domains with dense constellations of sources and observers in certain locations. For such problems using methods such as Fast Multiple Method \cite{GREENGARD1987325}, \cite{558669} or Non-Uniform Grid Interpolation method \cite{LI20108463} can be efficient. These methods can substitute the FFT-based approach for the near-zone component and the ideas presented here can be adapted to extend these methods to account for the additional images required due to periodicity.

\textcolor{black}{FFT-PIM can be used in a number of problems in electromagnetics, acoustics, and quantum mechanics. The evaluation of the periodic sums with dynamic (Helmholtz potential) PGF can be used in the context of electromanetic and acoustic integral equations for characterizing wave propagation, radiation, scattering, and dispersion diagrams for periodic arrays. The evaluation of the periodic sums with \textcolor{black}{NPSP} PGF can be used for modeling infinite periodic arrays or mimic infinite domains in electro-/magneto-statics, micromagnetics, and density functional theory. The evaluation of the periodic sums with static PGF with a phase shift can be used to calculate dispersion diagrams in micromagnetics, e.g., spin waves or linearized density functional theory.}   

\section{Code availability}
We open-sourced our far-zone component code package Periodic Unit Far Field (PUFF) under Apache-2.0 license on GitHub (\url{https://github.com/UCSD-CEM/PeriodicUnitFarField}).

\section{Acknowledgments}
This work was supported in part by the Quantum Materials for Energy Efficient Neuromorphic-Computing (Q-MEEN-C), an Energy Frontier Research Center funded by the U.S. Department of Energy, Office of Science, Basic Energy Sciences under Award No. DESC0019273. The work was also supported in part by Binational Science Foundation, grant \#2022346. The work used Purdue Anvil cluster at Rosen Center for Advanced Computing (RCAC) in Purdue University through allocation ASC200042 from the Advanced Cyberinfrastructure Coordination Ecosystem: Services \& Support (ACCESS) program \cite{10.1145/3569951.3597559}, which is supported by National Science Foundation grants \#2138259, \#2138286, \#2138307, \#2137603, and \#2138296.

\bibliographystyle{ieeetr}
\bibliography{references}

\end{document}